\documentclass[a4paper]{scrartcl}
\usepackage{amssymb}
\usepackage{amsmath}
\usepackage{amsthm}
\usepackage{latexsym}
\usepackage{enumerate}
\usepackage{graphicx,color}
\usepackage{hyperref}
\usepackage{dsfont}

\DeclareFontFamily{OT1}{pzc}{}
\DeclareFontShape{OT1}{pzc}{m}{it}{<-> s * [1.10] pzcmi7t}{}
\DeclareMathAlphabet{\mathpzc}{OT1}{pzc}{m}{it}

\usepackage{xcolor}
\usepackage{enumitem}

\usepackage[round,authoryear,longnamesfirst]{natbib}
\newcommand*{\doi}[1]{\href{http://dx.doi.org/#1}{doi: #1}}
\usepackage[autostyle=true]{csquotes} 

\renewcommand{\d}{\mathrm{d}}
\newcommand{\E}{\mathrm{E}}

\makeatletter
\newcommand*\bigcdot{\mathpalette\bigcdot@{.5}}
\newcommand*\bigcdot@[2]{\mathbin{\vcenter{\hbox{\scalebox{#2}{$\m@th#1\bullet$}}}}}
\makeatother

\numberwithin{equation}{section}

\numberwithin{figure}{section}

\theoremstyle{plain}
\newtheorem{theorem}{Theorem}[section]
\newtheorem{proposition}[theorem]{Proposition}
\newtheorem{corollary}[theorem]{Corollary}
\newtheorem{lemma}[theorem]{Lemma}

\theoremstyle{definition}
\newtheorem{definition}[theorem]{Definition}
\newtheorem{condition}[theorem]{Condition}

\newtheorem{remark}[theorem]{Remark}

\usepackage{tikz}
\usetikzlibrary{arrows,positioning,calc}
\usetikzlibrary{decorations.pathreplacing}
\tikzset{
    >=stealth',
    punkt/.style={
           rectangle,
           rounded corners,
           draw=black, thick,
           text width=5.5em,
           minimum height=2em,
           text centered},
    punktl/.style={
           rectangle,
           rounded corners,
           draw=black, thick,
           text width=7em,
           minimum height=2em,
           text centered},
    pil/.style={
           ->,
           shorten <=4pt,
       shorten >=4pt
    },
    pildotted/.style={
           ->,
           shorten <=4pt,
           shorten >=4pt,
  dotted,
  },
    punktf/.style={
           rectangle,
           text width=4.0em,
           minimum height=1.5em,
           text centered},
    punktfTop/.style={
           rectangle,
           text width=4.0em,
           minimum height=1.5em,
           text centered,
           append after command={
               [thick,shorten >=0.2bp, shorten <=0.2bp]
               (\tikzlastnode.north west)edge(\tikzlastnode.north east)
}
    },
    punktfBot/.style={
           rectangle,
           text width=4.0em,
           minimum height=1.5em,
           text centered,
           append after command={
               [thick,shorten >=0.2bp, shorten <=0.2bp]
               (\tikzlastnode.south west)edge(\tikzlastnode.south east)
            }
    }
}

\newcommand\xqed[1]{%
  \leavevmode\unskip\penalty9999 \hbox{}\nobreak\hfill
  \quad\hbox{#1}}
\newcommand\demormk{\xqed{$\triangledown$}}

\newcommand\democd{\xqed{$\diamond$}}
\newcommand\demodef{\xqed{$\triangle$}}

\usepackage{authblk}

\title{Dynamics of state-wise prospective reserves in the presence of non-monotone information}

\author[1]{Marcus C.~Christiansen}
\author[2,3,$\star$]{Christian Furrer}

\affil[1]{\footnotesize Institut f{\"u}r Mathematik, Carl von Ossietzky Universit{\"a}t Oldenburg,  Carl-von-Ossietzky-Stra{\ss}e 9--11, DE-26129 Oldenburg, Germany.}
\affil[2]{\footnotesize PFA Pension, Sundkrogsgade 4, DK-2100 Copenhagen \O, Denmark.}
\affil[3]{\footnotesize Department of Mathematical Sciences, University of Copenhagen, Universitetsparken 5, DK-2100 Copenhagen \O, Denmark.}
\affil[$\star$]{\footnotesize Corresponding author. E-mail: \href{mailto:furrer@math.ku.dk}{furrer@math.ku.dk}.}

\date{\vspace{-8mm}}

\begin{document}

\maketitle


\begin{abstract}

In the presence of monotone information, the stochastic Thiele equation describing the dynamics of state-wise prospective reserves is closely related to the classic martingale representation theorem. When the information utilized by the insurer is non-monotone, the classic martingale theory does not apply. By taking an infinitesimal approach, we derive a generalized stochastic Thiele equation that allows for information discarding. En passant, we solve some open problems for the classic case of monotone information. The results and their implication in practice are illustrated via examples where information is discarded upon and after stochastic retirement.

\end{abstract}

\vspace{5mm}


\vspace{5mm}

%


\vspace{5mm}

\section{Introduction}

Life insurers frequently employ reduced information in the valuation of liabilities due to e.g.\ legal constraints. The possibility of information discarding leads to potentially decreasing flows of information for which classic martingale-based methods do not apply. Based on the novel \textit{infinitesimal approach} proposed and developed in~\citet{Christiansen2020}, we study the dynamics of so-called state-wise prospective reserves in the presence of non-monotone information. Our main contribution is a generalization of the stochastic Thiele equation from~\citet{Norberg1992,Norberg1996} to allow for non-monotone information. Even in the presence of monotone information, \citet{Norberg1992,Norberg1996} leaves open some important technical questions regarding the concept of state-wise prospective reserves; an important secondary contribution of this paper is to answer such questions by constructing suitably regular modifications of the involved processes.

Non-monotone information structures appear in various contexts, including
\begin{itemize}
\item Legal constraints: By exercising the `right to be forgotten' according to Article 17 of the General Data Protection Regulation of the European Union (GDPR, \citeyear{GDPR}), the policyholder may ask the insurer to delete parts of e.g.\ health related data at discretion.
\item Risk-sharing: Markov chain models remain one of the most popular classes of multi-state life insurance models. If utilized, the required information becomes non-monotone. Practitioners might employ Markov chain models to effectively generate risk-sharing or pooling mechanisms which discriminate only on the basis of the current states of the insured. The validity of such practices might be analyzed by taking as given the non-monotone `as if Markovian' information structure rather than the probabilistic model (a Markov chain).
\item Notions of surplus: In \citet{Norberg1999} different notions of surplus based on non-monotone information are proposed, including so-called state-wise surplus. 
\item Inclusion of big data: Data from activity trackers, social media, etc.\ can improve forecasts of the mortality and morbidity of individual insured. Data providers might implement self-imposed information restrictions, such as the option of auto-deletion after certain time limits, to improve data privacy.
\end{itemize}
The stochastic Thiele equation describing the dynamics of state-wise prospective reserves is central to the mathematics of multi-state life insurance. Besides providing a convenient starting point for elegant derivations of Feynman-Kac formulas, leading to
\begin{itemize}
\item efficient computational schemes for the calculation of prospective reserves,
\end{itemize}
the study of dynamics of (state-wise) prospective reserves is crucial in relation to e.g.
\begin{itemize}
\item Cantelli's theorem and the definition of prospective reserves in the presence of reserve-dependent payments, see e.g.~\citet{MilbrodtStracke1997} and \citet{ChristiansenDjehiche2020}.
\item Hattendorff's theorem on the (lack of) correlation between losses, see e.g.~\citet{RamlauHansen1988a} and \citet{Norberg1992}.
\item The emergence and decomposition of surplus and safe-side calculations and sensitivity analyses applied in the context of risk margins and the construction of technical bases, see e.g.~\citet{RamlauHansen1988b}, \citet{Linnemann1993}, \citet{Norberg2001}, \citet{KalashnikovNorberg2003}, \citet{Christiansen2008a,Christiansen2008b}, and \citet{ChristiansenSteffensen2013}.
\end{itemize}
Besides the generalization of Hattendorff's theorem found in~\citet{Norberg1992} and the contributions of~\citet{ChristiansenDjehiche2020}, the aforementioned concepts and methods have primarily been developed for Markov chain models and semi-Markov models. By revisiting and improving upon the work in~\citet{Norberg1992,Norberg1996} on the stochastic Thiele equation in the presence of monotone information, we take an important step towards closing technical and conceptual gaps in the mathematics of multi-state life insurance. By further allowing for information discarding following~\citet{Christiansen2020}, we also seek to initialize a research program that aims to advance the mathematics of multi-state life insurance above and beyond martingale-based methods.

The multi-state approach to life insurance dates back at least to \citet{Hoem1969}, where Thiele equations describing the dynamics of the state-wise prospective reserves are derived for Markov chain models. These differential equations were revisited by Norberg in his seminal paper of 1991~\citep{Norberg1991} and have been generalized in various directions. This includes relaxing the assumption of Markovianity to allow for duration dependence (semi-Markovianity), taking market risks into account, and the study of higher order moments of prospective reserves, see e.g.~\citet{Moller1993}, \citet{Steffensen2000}, \citet{Helwich2008}, \citet{AdeChristiansen2017}, and \citet{Bladt2020}.  We should mention that
while the approach in~\citet{Steffensen2000} is very general, the results are only established under strict smoothness conditions that might not be satisfied in practice. In this paper, and contrary to~\citet{Steffensen2000}, the only source of randomness consists of the state of the insured, which is modeled as a non-explosive pure jump process on a countable state space.

The ordinary Thiele equations are essentially Feynman-Kac type results. In contrast, the stochastic Thiele equation from~\citet{Norberg1992,Norberg1996} is a stochastic differential equation that applies irregardless of the intertemporal dependence structure and reveals the universality of Thiele's original equation. In the presence of monotone information, the dynamics of prospective reserves are characterized by identifying integrands in the classic martingale representation theorem for random counting measures \citep{Norberg1992,Norberg1996,ChristiansenDjehiche2020}. In similar fashion, our approach relies on the infinitesimal martingale representation theorem from~\citet{Christiansen2020}, which extends the classic martingale representation theorem for random counting measures to allow for non-monotone information. Furthermore, we develop for the first time a mathematically sound concept of state-wise prospective reserves which applies irregardless of the intertemporal dependence structure.

The paper is structured as follows. In Section~\ref{Section:sketch}, a short sketch of concepts and results is given. In Section~\ref{sec:setup}, we present the probabilistic setup, and in Section~\ref{sec:setup_stoc_re}, this setup is exemplified by considering information discarding upon and after stochastic retirement. Mathematically sound concepts of state-wise prospective reserves are developed in Section~\ref{sec:state_wise_prores}, while Section~\ref{sec:inf_comps} discusses the concept of infinitesimal compensators within the present setup. Piecing everything together by an application of the infinitesimal martingale representation theorem, Section~\ref{sec:stoch_thiele} contains the derivation of the main result of this paper: the generalization of the stochastic Thiele equation to allow for non-monotone information. We conclude by returning to examples of information discarding upon and after stochastic retirement, where we study the general dynamics of prospective reserves (Section~\ref{sec:examples}) and derive Feynman-Kac formulas beyond the classic cases of Markovianity and semi-Markovianity (Section~\ref{sec:examples2}).

\section{Short sketch of concepts and results}\label{Section:sketch}

This section quickly makes the central ideas and results of this paper comprehensible by looking at a simplified model setting, by avoiding technicalities, and by postponing the proofs. In particular, we are not explicit about completion of $\sigma$-algebras, we omit the inclusion of point probability mass in the distributions, and we require sojourn and transition payments to be deterministic. The subsequent sections are then aimed at making the mathematics rigorous while also covering a wider range of model settings.

Suppose that the random pattern of states of an individual insured is modeled by a pure jump process $Z=(Z_t)_{t\geq 0}$ with values in a countable state space $S$. The total information available is given by the filtration $\mathcal{F}=(\mathcal{F}_t)_{t\geq0}$ naturally generated by $Z$. We relate to the random pattern of states $Z$ the counting processes
$N_{jk}(t)=(N_{jk}(t))_{t\geq0}$, $j,k \in S$, $j \neq k$, giving the number of jumps of $Z$ from state $j$ to state $k$:
\begin{align}\label{DefNjk}
N_{jk}(t) :=  \#\left\{
s \in (0,t] : Z_{s-} = j, Z_s = k
\right\}\!,\hspace{10mm} t\geq 0.
\end{align}
Let $B=(B(t))_{t \geq 0}$ represent the accumulated benefits minus premiums of a life insurance contract.
Suppose that
\begin{align*}
\d B( t)&=  \sum_{j\in S} \mathbf{1}_{\{Z_{t-}=j\}} b_j(t) \, \d t
+
\sum_{j,k\in S\atop j \neq k} b_{jk}(t) \, \d N_{jk}( t),\quad B(0)=0,
\end{align*}
where $b_j$ and $b_{jk}$ are bounded functions that describe sojourn and transition payments.
The suitably discounted accumulated future payments $Y=(Y(t))_{t\geq0}$ are given by
\begin{align}\label{DefOfY}
Y(t)=
\int_{(t,\infty)} e^{-\int_t^s r(u)\, \d u}\,  \d B( s), \quad t\geq0,
\end{align}
where $r$ is a deterministic short rate, and the prospective reserve at time $t$ is defined as
\begin{align}\label{PRE:DefProspRes0}
V^{\mathcal{F}}(t) = \E[ Y(t) \, | \,  \mathcal{F}_{t}], \quad t\geq0.
\end{align}
Equation~\eqref{PRE:DefProspRes0} defines the prospective reserve only point-wise almost surely, which means that the prospective reserve is not necessarily well-defined as a stochastic process.  Yet, $V^{\mathcal{F}}$ has a unique c\`{a}dl\`{a}g modification known as the optional projection for which the Ito differential $\d V^{\mathcal{F}}$ exists.

\subsubsection*{Full information setting}

In the classic literature on multi-state life insurance, $Z$ is typically assumed to be a Markov process such that $V^{\mathcal{F}}(t)=\E[ Y(t) \,  | \,  Z_t]$, and state-wise prospective reserves are then defined as $V_j(t)=\E[ Y(t) \,  | \,  Z_t=j]$, $j \in S$. Without the Markov assumption, the general definition of state-wise prospective reserves is largely unclear. In~\citet{Norberg1996}, a definition of state-wise prospective reserves is suggested in case the information $\mathcal{F}$ can be replaced by the current information of an expanded state process. More precisely, it is assumed that there exist processes $F^j=(F^j_t)_{t\geq 0}$, $j \in S$, with values in some appropriate state space $E$ such that
the bivariate process $(Z_t,F_{t})_{t \geq 0}$ for $F_t:= \sum_{j\in S} \mathbf{1}_{\{Z_t=j\}} F^j_t$ satisfies
\begin{align*}
V^{\mathcal{F}}(t)=\E[ Y(t) \, | \,  Z_t,F_{t}],  \quad t\geq 0,
\end{align*}
and state-wise prospective reserves are then defined according to
\begin{align}\label{PRE:DefProspRes}
  V^{\mathcal{F}}_j(t) = \E[ Y(t) \,  | \,  Z_t=j,F^j_{t}], \quad j\in S, t\geq 0,
\end{align}
since this definition implies  that
\begin{align}\label{DefiningEqualityVFj}
V^{\mathcal{F}}(t)= \sum_{j\in S} \mathbf{1}_{\{Z_t=j\}} V^{\mathcal{F}}_j(t), \quad t\geq0.
\end{align}
The conditional expectation in \eqref{PRE:DefProspRes} shall be read as the factorized conditional expectation $\E[ Y(t) \, | \,  (Z_t, F^j_{t})= \cdot \, ]$  evaluated at $(j,F^j_t)$. We use such short-hand notation throughout this paper.
It is left open in \citet{Norberg1996} whether a suitable set of processes $F^j$, $j \in S$, always exists. This paper presents a construction for $F^j$, $j\in S$, such that we even have
\begin{align}\label{PRE:FRepresentation}
  \mathcal{F}_t = \sigma( Z_t,F_{t}), \quad t \geq 0,
\end{align}
which implies that definition \eqref{PRE:DefProspRes} works irregardless of the intertemporal dependence structure of $Z$.
Unfortunately, since \eqref{PRE:DefProspRes} defines $V^{\mathcal{F}}_j$, $j\in S$, only point-wise almost surely, the state-wise prospective reserves are not necessarily well-defined as stochastic processes, and, moreover, it is unclear whether the Ito differentials $\d  V^{\mathcal{F}}_j$, $j\in S$, exist. Note that the Ito differentials are needed in the stochastic Thiele equation. In~\citet{Norberg1996}, the questions concerning well-definedness and existence are left open. In this paper, we construct modifications of $V^{\mathcal{F}}_j$, $j\in S$, such that
$t \mapsto \mathbf{1}_{\{Z_{t-}=j\}}\d V^{\mathcal{F}}_j(t)$, $j \in S$, exist. As a consequence hereof, we can give  the first mathematically rigorous proof of the stochastic Thiele equation in a general non-Markovian setting:
\begin{align}\label{PRE:Thiele}
      0=\sum_{j\in S} \mathbf{1}_{\{Z_{t-}=j\}} \bigg(\d  V^{\mathcal{F}}_j( t) -  V^{\mathcal{F}}_j(t)\, r(t)\, \d t  + b_j(t)\, \d t+  \sum_{k:k\neq j}  R^{\mathcal{F}}_{jk}(t) \,\mu^{\mathcal{F}}_{jk}(t)\, \d t\bigg).
\end{align}
In~\citet{Norberg1992}, 
the equation~\eqref{PRE:Thiele} is denoted as a \textit{stochastic Thiele's differential equation}. Under appropriate smoothness conditions,  the transition intensity $\mu^{\mathcal{F}}_{jk}(t)$ equals
\begin{align}\label{PRE:DefMu}
\mu^{\mathcal{F}}_{jk}(t) = \lim_{h\downarrow 0}\frac{1}{h} P(Z_{t+h}=k \,|\, Z_t=j,F^j_t).
\end{align}
The so-called sum at risk
\begin{align}\label{PRE:SumAtRisk} \begin{split}
  R^{\mathcal{F}}_{jk}(t)= b_{jk}(t)&+ \E\big[  Y(t) \, \big| \, Z_{t-}=j,F^j_{t-},Z_t=k \big] - \E\big[ Y(t) \, \big| \, Z_{t-}=j,F^j_{t-},Z_t=j\big]
\end{split}\end{align}
describes the change in the reserve upon a jump of $Z$ from $j$ to $k$ at time $t$.  If $Z$ is Markovian, then we obtain $R^{\mathcal{F}}_{jk}(t)=b_{jk}(t) + V^{\mathcal{F}}_k(t) -V^{\mathcal{F}}_j(t)$. In case $Z$ is not Markovian, we still get
\begin{align}\label{PRE:SumAtRiskSimplified}
 \mathbf{1}_{\{Z_{t-}=j\}} \, R^{\mathcal{F}}_{jk}(t)  = \mathbf{1}_{\{Z_{t-}=j\}} \Big(b_{jk}(t) + V^{\mathcal{F}}_k(t) -V^{\mathcal{F}}_j(t)\Big) 
\end{align}
if the state-wise prospective reserves are carefully defined, as we show in this paper.

\subsubsection*{Non-monotone information setting}
In order to model information discarding, in a next step we introduce processes $G^j$, $j \in S$, so that
\begin{align*}
  \mathcal{G}_t = \sigma( Z_t,G_t) \subseteq \mathcal{F}_t, \quad t \geq 0,
\end{align*}
for $G_t := \sum_{j\in S} \mathbf{1}_{\{Z_t=j\}} G^j_t$, where the sequence of sub-$\sigma$-algebras $\mathcal{G}=(\mathcal{G}_t)_{t \geq 0}$ describes the reduced information. While $\mathcal{F}$ is a filtration, $\mathcal{G}$ is non-monotone unless $\mathcal{G}=\mathcal{F}$. We aim to study the dynamics of the prospective reserve under information $\mathcal{G}$, which is defined as
\begin{align}\label{eq:classic_res_def}
V(t) = \E[ Y(t) \,  | \,  \mathcal{G}_{t}], \quad t\geq0.
\end{align}
Equation~\eqref{eq:classic_res_def} defines the prospective reserve under information $\mathcal{G}$ only point-wise almost surely, yet even under non-monotone information it has a unique c\`{a}dl\`{a}g modification for which the Ito differential $\d V$ exists, cf.~Theorem 4.1 in~\citet{Christiansen2020}. The analogue of~\eqref{PRE:DefProspRes} for $\mathcal{G}$ instead of $\mathcal{F}$ is
\begin{align}\label{PRE:DefProspResG}
  V_j(t) = \E[ Y(t) \, | \,  Z_t=j , G^j_{t}], \quad j\in S, t\geq0,
\end{align}
which yields the equation
\begin{align}\label{DefiningPropertyVj}
V(t)= \sum_{j\in S} \mathbf{1}_{\{Z_t=j\}} V_j(t), \quad t\geq0.
\end{align}
Again we have the problem that \eqref{PRE:DefProspResG} defines $V_j$, $j\in S$, only point-wise almost surely. This paper constructs modifications of $V_j$, $j \in S$, such that $t \mapsto \mathbf{1}_{\{Z_{t-}=j\}}\d V _j(t)$, $j \in S$, exist.  Based on that construction, we then discover a generalized stochastic Thiele equation for the non-monotone information setting:
\begin{align}\label{PRE_generalizedThiele}\begin{split}
      0=\sum_{j\in S} \mathbf{1}_{\{Z_{t-}=j\}} \bigg(\d  V _j( t) -  V _j(t)\, r(t)\, \d t  + b_j(t)\, \d t&+  \sum_{k:k\neq j} \int_E R _{jk}(t,g) \,\mu _{jk}(t,\d g)\, \d t\\
      &+  \sum_{k:k\neq j} \int_E \overline{R} _{kj}(t,g) \,\overline{\mu} _{kj}(t,\d g)\, \d t\bigg),
\end{split}\end{align} 
where, under appropriate smoothness conditions,
\begin{align*}
  \mu _{jk}(t, A) = \lim_{h\downarrow 0}\frac{1}{h} P\big(Z_{t+h}=k, G^k_{t+h}\in  A \, \big| \, Z_t=j,G^j_t\big)
\end{align*}
describes the transition intensity for a jump from $Z_t=j$ to $Z_t=k$ with $G^k_{t}\in  A$ at time $t$,
\begin{align*}
  R _{jk}(t,g)= b_{jk}(t) &+ \E\big[   Y(t) \, \big| \, Z_{t-}=j,G^j_{t-},Z_t=k,G^k_{t}=g \big] - \E\big[ Y(t) \, \big| \, Z_{t-}=j,G^j_{t-},Z_t=j\big]
\end{align*}
describes the sum at risk for a transition from $Z_{t-}=j$ to $Z_t=k$ with $G^k_{t}=g$ at time $t$,
\begin{align*}
  \overline{\mu} _{kj}(t, A) = \lim_{h\downarrow 0}\frac{1}{h} P\big(Z_{t-h}=k,G^k_{t-h}\in  A \,  \big| \, Z_t=j,G^j_t\big)
\end{align*}
describes a backward-looking transition intensity for the backward step from $Z_t=j$ to $Z_{t-}=k$ with $G^k_{t-}\in A$ at time $t$, and
\begin{align*}
  \overline{R} _{kj}(t,g)= &\,\E\big[ Y(t) \, \big| \, Z_{t-}=k,G^k_{t-}=g,Z_t=j,G^j_{t} \big]  - \E\big[ Y(t) \, \big| \, Z_{t-}=j,Z_t=j,G^j_{t}\big]
\end{align*}
describes a backward-looking sum at risk for the backward step from $Z_t=j$ to $Z_{t-}=k$ with $G^k_{t-}=g$ at time $t$.

The second line in \eqref{PRE_generalizedThiele} describes the effect of information discarding on the prospective reserve. It vanishes in particular if $\mathcal{G}$ is monotone, i.e.~when $\mathcal{G}=\mathcal{F}$. The integral over $E$ in the first line of \eqref{PRE_generalizedThiele} is required since our model allows us to add not only the information $\sigma(Z_t)$ in the step from $\mathcal{G}_{t-}$ to $\mathcal{G}_t$ but also further information from $\mathcal{F}_{t-}$ not already included in $\mathcal{G}_{t-}$. This scenario occurs if previously discarded information is recovered at time $t$. In case of
$ \mathcal{G}_{t}\subseteq  \mathcal{G}_{t-} \vee \sigma(Z_t)$ for all $t$,
the first line of \eqref{PRE_generalizedThiele} takes the same form as \eqref{PRE:Thiele} with $\mu _{jk}$ and $R _{jk}$ defined analogously to \eqref{PRE:DefMu} and \eqref{PRE:SumAtRisk} but with $G^j$ in place of $F^j$.

\subsubsection*{Example: Stochastic retirement}

We consider a multi-state life insurance contract where the insurer distinguishes between two different health statuses of an insured, here modeled as the states $\{1, 2\}$. At a random time $\eta$ the observation of the health status is stopped and the past health records are discarded, here modeled as a jump to state $\text{p}$ at time $\eta$. When the insured dies, here given as a jump to an absorbing state $\text{d}$, past health records are also discarded. In both cases, the time of retirement $\eta$ is not recorded. So our state space is $S=\{1,2,\text{p},\text{d}\}$, the processes $G^1,G^2$ equal $F^1,  F^2$, and
the processes $G^{\text{p}}$ and  $G^{\text{d}}$ are redundant. In this paper, we show that it is often unreasonable to assume $Z$ to be a Markov process. In particular, $\E[Y(t) \, | \, \mathcal{F}_t]$ is indeed different from $\E[Y(t) \, | \,Z_t]$.

Our model setting occurs e.g.\ in the context of pensions with disability coverage. Before retirement, the insurer distinguishes between active and disabled insured, but upon and after retirement, all insured are just viewed as pensioners without regard to their current health status or past health records. The time point $\eta$ of information discarding is indeed random if the retirement age is flexible. It is common in actuarial practice to compute the state-wise prospective reserve in state $\text{p}$ without the health record $F^{\text{p}}$, but then the classic Thiele differential equations do not apply anymore, as we now show.

If $Z$ was to be a Markov process, then $V^{\mathcal{F}}_j(t)$ and $\mu^{\mathcal{F}}_{jk}(t)$ are deterministic, and the classic Thiele differential equations for $V^{\mathcal{F}}_j$, $j \in S$, can be directly derived from~\eqref{PRE:Thiele}. If $Z$ is not Markovian, then further steps are needed in order to transform~\eqref{PRE:Thiele} into a deterministic system.  For the sake of simplicity, we assume here the Markovianity of an extended state process $\tilde{Z}$, defined as the process where state $\text{p}$ is split into two sub-states $\{3,4\}$ giving the (unobservable) health statuses of the insured upon and after retirement. One can then show that the processes $V^{\mathcal{F}}_{\text{p}}$ and $\mu^{\mathcal{F}}_{\text{p}\text{d}}$ are of the form
\begin{align*}
  V^{\mathcal{F}}_{\text{p}}(t) &= V^{\mathcal{F}}_{\text{p}}(t, \eta, Z_{\eta-}),\\
  \mu^{\mathcal{F}}_{\text{p}\text{d}}(t) &= \mu^{\mathcal{F}}_{\text{p}\text{d}}(t, \eta, Z_{\eta-})
\end{align*}
for  deterministic functions $V^{\mathcal{F}}_{\text{p}}(t, r,\ell)$  and  $\mu^{\mathcal{F}}_{\text{p}\text{d}}(t, r, \ell)$ with $r \in [0,t]$ and $\ell \in \{1,2\}$. If we disregard the past health status and the time of retirement for an insured in state $\text{p}$, then we are basically replacing $V^{\mathcal{F}}_{\text{p}}(t, \eta ,Z_{\eta-})$  and  $\mu^{\mathcal{F}}_{\text{p}\text{d}}(t, \eta, Z_{\eta-})$ by the $\mathcal{G}$-averages
\begin{align*}
 V_{\text{p}}(t)&=\E[ V^{\mathcal{F}}_{\text{p}}(t, \eta ,Z_{\eta-}) \, | \, Z_t=\text{p}]   ,\\
 \mu _{\text{p}\text{d}}(t)&= \E[ \mu^{\mathcal{F}}_{\text{p}\text{d}}(t,\eta, Z_{\eta-}) \, | \, Z_{t-}=\text{p}].
\end{align*}
All other state-wise prospective reserves are deterministic and satisfy the equation
  \begin{align*}
  V_j(t) = V^{\mathcal{F}}_j(t), \quad j \in \{1,2,\text{d}\}, t\geq 0.
\end{align*}
By applying~\eqref{PRE:Thiele}, we are able to derive the following classic Thiele differential equations:
\begin{align*}
      \d  V^{\mathcal{F}}_j( t) &=  Y^{\mathcal{F}}_j(t)\, r(t)\, \d t  - b_j(t)\, \d t- \hspace{-1mm}   \sum_{k \in \{1,2,\text{d}\} \atop j \neq k} R^{\mathcal{F}}_{jk}(t) \,\mu^{\mathcal{F}}_{jk}(t)\, \d t
      - R^{\mathcal{F}}_{j\text{p}}(t,t,j)\,\mu^{\mathcal{F}}_{j\text{p}}(t), \quad j \in \{1,2,\text{d}\},\\
      \d  V^{\mathcal{F}}_{\text{p}}( t,r,\ell) &=  V^{\mathcal{F}}_{\text{p}}(t,r,\ell)\, r(t)\, \d t  - b_{\text{p}}(t)\, \d t-  R^{\mathcal{F}}_{\text{p}\text{d}}(t,r,\ell) \,\mu^{\mathcal{F}}_{\text{p}\text{d}}(t,r,\ell)\, \d t, \quad r \leq t,\,\ell\in \{1,2\},
\end{align*}
where the sum at risks are given by
\begin{align*}
R^{\mathcal{F}}_{jk}(t)&= b_{jk}(t) + V^{\mathcal{F}}_{k}( t)-V^{\mathcal{F}}_{j}( t), \quad j,k \in \{1,2,\text{d}\}, j \neq k,\\
  R^{\mathcal{F}}_{j\text{p}}(t,t,j)&= b_{j\text{p}}(t) + V^{\mathcal{F}}_{\text{p}}( t,t,j)-V^{\mathcal{F}}_{j}( t),\quad j \in \{1,2\},\\
  R^{\mathcal{F}}_{\text{p}\text{d}}(t,r,\ell) &= b_{\text{p}\text{d}}(t) + V^{\mathcal{F}}_{\text{d}}( t)-V^{\mathcal{F}}_{\text{p}}( t,r,\ell).
\end{align*}
By applying \eqref{PRE_generalizedThiele}, we can derive for $V_{\text{p}}$ the deterministic differential equation
\begin{align}\label{PRE:VG3equation}\begin{split}
  \d  V _{\text{p}}(t) =& \, V _{\text{p}}(t)\, r(t)\,\d t- b_{\text{p}}(t) \, \d t- \big( b_{\text{p}\text{d}}(t)  + V_{\text{d}}(t) - V_{\text{p}}(t)\big) \,\mu _{\text{p}\text{d}}(t)\,\d t\\
   &+ \big( V^{\mathcal{F}}_{\text{p}}(t,t,1) - V_{\text{p}}(t)\big) \, \overline{\mu} _{1{\text{p}}}(t) \, \d t+ \big( V^{\mathcal{F}}_{\text{p}}(t,t,2) - V_{\text{p}}(t)\big)  \,\overline{\mu} _{2{\text{p}}}(t)\, \d t.
\end{split}\end{align}
The second line in \eqref{PRE:VG3equation} relates to information discarding. It consists of backward looking sum at risks and the backward looking transition intensities
 \begin{align}\label{PRE:OverlineLambda}
   \overline{\mu} _{j\text{p}}(t) = \lim_{h\downarrow 0}\frac{1}{h} P(Z_{t-h}=j \, | \, Z_t=\text{p} )= \frac{P(Z_t=j)}{P(Z_t=\text{p})}\mu^{\mathcal{F}}_{j\text{p}}(t), \quad j \in \{1,2\}.
 \end{align}
   If $Z$ is a  Markov process, then the second line in \eqref{PRE:VG3equation} vanishes and $\mu_{\text{p}\text{d}}^{\mathcal{F}}(t,r,\ell) \equiv \mu_{\text{p}\text{d}}(t)$, which leads us back to the classic Thiele differential equations.


\section{Modeling of non-monotone information}\label{sec:setup}



For a multitude of reasons, including legal constrains, the inclusion of big data, or to achieve model simplifications, the insurer might not have access to or desire to utilize all information available to it. Examples include the newly introduced General Data Protection Regulation of the European Union (GDPR, \citeyear{GDPR}), where Article 17 describes a so-called `right to be forgotten', and the restriction to an `as if Markovian' information structure even when the Markov property is not satisfied. Representing the resulting utilized information as a sub-sequence of $\sigma$-algebras, one typically finds that the sequence is non-monotone because certain pieces of information are discarded en route.

In this section, we introduce a general modeling framework for multi-state life insurance settings that involve information discarding. The framework is strongly related to the general theory of non-monotone information for jump processes introduced in~\citet{Christiansen2020}. In the following section, an application regarding stochastic retirement is discussed. This leads to the specification of some explicit cases of non-monotone information that serve as the main examples in the ensuing investigation.

Let $(\Omega,\mathbb{F},P)$ be a complete probability space with null-sets $\mathcal{A}$, and let $Z=(Z_t)_{t\geq 0}$ be a random pattern of states (pure jump process) on a countable state space $S$ giving at each time $t$ the state of the insured in $S$. Let
$$I_j(t) := \mathbf{1}_{\{Z_t=j\}}, \quad j \in S, t\geq0,$$
 be the corresponding state indicator processes. Denote by $N_{jk}=(N_{jk}(t))_{t\geq0}$, $j,k\in S$, $j \neq k$, the counting processes associated with $Z$; they are rigorously defined in~\eqref{DefNjk}. We generally assume that
\begin{align}\label{assump:tech_noexplo}
  \E\bigg[\sum_{j,k\in S, j\neq k}\! N_{jk}(t) \,\bigg] < \infty, \quad t\geq0,
\end{align}
since this ensures that compensated counting processes are true martingales. Furthermore, the assumption is sufficient (but not necessary) to guarantee that $Z$ is non-explosive:
\begin{lemma}\label{FinitelyManyJumps}
For each $t \geq 0$ and almost each $\omega \in \Omega$, the process $Z$ has at most finitely many jumps on $[0,t]$.
\end{lemma}
\begin{proof}
If $Z$ had infinitely many jumps on $[0,t]$ with positive probability, then the expectation in \eqref{assump:tech_noexplo} would be infinite
\end{proof}

In the following, stochastic processes are generally only defined and characterized up to evanescence.

\subsubsection*{Full information setting} 

The total information available is denoted $\mathcal{F}=(\mathcal{F}_t)_{t\geq 0}$; it is the right-continuous
filtration naturally generated by $Z$ added null-sets (to make it complete). Its components are given by
\begin{align*}
\mathcal{F}_t := \sigma(Z_s : s \leq t) \vee \mathcal{A}
\end{align*}
for $t\geq 0$. Since $\mathcal{F}$ is a filtration, it represents monotone (non-decreasing) information. We now construct a Markovian representation of $\mathcal{F}$, i.e.\ we define a jump process whose current state represents the full information $\mathcal{F}$.

Let $J = (J_t)_{t\geq0}$ be the counting process that gives the number of times that the process $Z$ has left the current state so far:
\begin{align*}
  J_t&:= \sum_{j \in S} I_j(t)\, J^j_t,\quad   J^j_t:=  \sum_{k : k \neq j} N_{jk}(t-) ,\quad j \in S,\, t \geq 0.
\end{align*}
Here $N_{jk}(0-) := 0$ for $j,k \in S$, $j \neq k$. The bivariate  jump process
$$\mathcal{Z}:=(Z,J)$$
 takes values in the state space
$\mathcal{S}:= S \times  \mathbb{N}_0$. Note that each state in $\mathcal{S}$ is visited by $\mathcal{Z}$ at most once. The entry and exit times are given by
\begin{align*}
T_x:= \inf\{t \geq 0 : \mathcal{Z}_t=x\}, \quad S_x:= \sup\{t \geq 0 : \mathcal{Z}_t=x\}, \quad x\in\mathcal{S},
\end{align*}
where we use the  conventions $\inf \emptyset := \infty$ and $\sup \emptyset := \infty$. We also define state indicator processes for the expanded jump process $\mathcal{Z}$ according to
\begin{align*}
\mathcal{I}_x(t) := \mathbf{1}_{\{\mathcal{Z}_t=x\}} = \mathbf{1}_{\{T_x \leq t < S_x\}}, \quad x\in\mathcal{S},t\geq 0,
\end{align*}
as well as counting processes via
\begin{align*}
\mathcal{N}_{xy}(t) :=  \#\left\{
s \in (0,t] : \mathcal{Z}_{s-} = x, \mathcal{Z}_s = y \right\}\!,\quad x,y\in\mathcal{S},x\neq y, t\geq 0.
\end{align*}
Set $\xi_x:=  (T_y\mathbf{1}_{\{T_y \leq T_x \}})_{y \in \mathcal{S}}$ for $x\in\mathcal{S}$. We interpret $\sigma(\xi_x)$ as the state-wise remaining information of $\mathcal{F}_t$ that is not already covered by $\sigma(\mathcal{Z}_t) $ upon $\mathcal{Z}_t=x$. This interpretation is a natural consequence of the following result.
\begin{lemma}\label{FTSZetaT}
For $x \in \mathcal{S}$ it holds that
\begin{align*}
\sigma(Z_s : s \leq t) \cap  \{ \mathcal{Z}_t=x\} &=  \sigma( \xi_x )  \cap  \{ \mathcal{Z}_t=x\}, \quad t\geq0, \\
\sigma(Z_s : s < t) \cap  \{ \mathcal{Z}_{t-}=x\} &=  \sigma( \xi_x )  \cap  \{ \mathcal{Z}_{t-}=x\}, \quad t>0.
\end{align*}
\end{lemma}
\begin{proof}
On $\{ \mathcal{Z}_t=x\} =\{T_x \leq t < S_x\}$ and $\{T_x < t \leq S_x\}$
the random vector $\xi_x$ provides the jump times and the departure states of $Z$ on $[0,t]$ and $[0,t)$, respectively.
\end{proof}
The $\sigma$-algebra $\sigma( \xi_x ) $ does not depend on time $t$. This technical feature will be crucial later on in order to construct modifications of the state-wise prospective reserves that have nice path properties.
For $\xi_{\mathcal{Z}_t}:= \sum_{x \in \mathcal{S}} \mathcal{I}_x(t) \,\xi_x$, $t\geq0$, Lemma~\ref{FTSZetaT} implies that
\begin{align}\label{FRepresentZXi}
\begin{split}
  \mathcal{F}_t &= \sigma (\mathcal{Z}_t,\xi_{\mathcal{Z}_t}) \vee \mathcal{A}, \quad t \geq 0,\\
  \mathcal{F}_{t-} &= \sigma (\mathcal{Z}_{t-},\xi_{\mathcal{Z}_{t-}}) \vee \mathcal{A}, \quad t > 0.
\end{split}\end{align}
This means that the process $(\mathcal{Z}_t,\xi_{\mathcal{Z}_t})_{t \geq 0}$ provides a Markovian representation of the full information $\mathcal{F}$. In a next step, we construct a Markovian representation of $\mathcal{F}$ that bases on $Z$ instead of $\mathcal{Z}$. To this end, let
\begin{align*}
F^j_t :=  \sum_{n \in \mathbb{N}_0} \mathbf{1}_{\{J^j_t=n\}}\, (n, \xi_{(j,n)}), \quad j\in S, t\geq 0.
\end{align*}
We interpret -- in accordance with the following result -- the information $\sigma(F_t^j)$ as the state-wise remaining information of $\mathcal{F}_t$ that is not already covered by $\sigma(Z_t)$ upon $Z_t = j$.
\begin{lemma}\label{FZF}
For $j \in S$ it holds that
\begin{align*}
\sigma(Z_s : s \leq t) \cap \{Z_t=j\} &= \sigma(F^j_t ) \cap \{Z_t=j\}, \quad t\geq0, \\
\sigma(Z_s : s < t) \cap \{Z_{t-}=j\} &= \sigma(F^j_{t-} ) \cap \{Z_{t-}=j\}, \quad t>0.
\end{align*}
\end{lemma}
\begin{proof}
Since $\{Z_t=j\}\cap \{J^j_t=n\} = \{T_{(j,n)} \leq t < S_{(j,n)}\}$ and $\{Z_{t-}=j\}\cap\{J^j_{t-}=n\} =\{T_{(j,n)} < t \leq S_{(j,n)}\}$, the statement follows from Lemma~\ref{FTSZetaT}.
\end{proof}
For $F_t := \sum_{j\in S} I_j(t) \, F_t^j$, $t\geq0$, Lemma~\ref{FZF} implies the following Markovian representation of the full information $\mathcal{F}$:
\begin{align*}
  \mathcal{F}_t &= \sigma(Z_t,F_t) \vee \mathcal{A}, \quad t \geq 0,\\
  \mathcal{F}_{t-} &= \sigma(Z_{t-},F_{t-}) \vee \mathcal{A}, \quad t > 0.
\end{align*}

\subsubsection*{Non-monotone information setting} \label{subsec:nonmono}

We model non-monotone information by constructing sub-$\sigma$-algebras $\mathcal{G}_t \subseteq \mathcal{F}_t$, $t\geq0$, as follows. Let  $ \zeta_x:(\Omega, \mathbb{F}) \rightarrow (E',\mathfrak{B}(E'))$, $x\in\mathcal{S}$, be random variables with values in a separable and completely metrizable topological space $E'$ (where $\mathfrak{B}(E')$ denotes the corresponding Borel $\sigma$-algebra) such that
\begin{align*}
 \sigma( \zeta_x) \subseteq  \sigma( \xi_x), \quad x \in \mathcal{S}.
\end{align*}
In the same way that $\xi_x$ gives the full information at time $t$ upon $\mathcal{Z}_t=x$, the quantity  $\zeta_x$ shall represent the reduced information at time $t$ upon $\mathcal{Z}_t=x$. Note that this implies that information is discarded only at the jump times of $Z$.

Analogously to~\eqref{FRepresentZXi}, we set
\begin{align}\label{DefG}\begin{split}
      \mathcal{G}_t &:=  \sigma( \mathcal{Z}_t,\zeta_{\mathcal{Z}_t})  \vee \mathcal{A},\quad t \geq 0,\\
       \mathcal{G}_{t-} &:= \sigma( \mathcal{Z}_{t-},\zeta_{\mathcal{Z}_{t-}})  \vee \mathcal{A}, \quad t>0.
\end{split}\end{align}
Our construction implies that the process $(\mathcal{Z}_t,\zeta_{\mathcal{Z}_t})_{t \geq 0}$ provides a Markovian representation of the reduced information $\mathcal{G}$. Similarly to Lemma~\ref{FTSZetaT}, we obtain
\begin{align}\begin{split}\label{GcondSigmaZeta}
 \sigma( \mathcal{Z}_t,\zeta_{\mathcal{Z}_t}) \cap   \{ \mathcal{Z}_t=x\} &=  \sigma( \zeta_x )  \cap  \{ \mathcal{Z}_t=x\}, \quad t \geq 0, \\
  \sigma( \mathcal{Z}_{t-},\zeta_{\mathcal{Z}_{t-}}) \cap   \{ \mathcal{Z}_{t-}=x\} &=  \sigma( \zeta_x )  \cap  \{ \mathcal{Z}_{t-}=x\}, \quad t > 0.
  \end{split}
\end{align}
In a next step, we construct a Markovian representation of $\mathcal{G}$ that bases on $Z$ instead of $\mathcal{Z}$. To this end, let
\begin{align*}
G^j_t :=  \sum_{n \in \mathbb{N}_0} \mathbf{1}_{\{J^j_t=n\}}\, (n, \zeta_{(j,n)}), \quad j\in S, t \geq 0.
\end{align*}
Note that this defines jump processes on the state space $E:=\mathbb{N}_0 \times E'$. Once again we might interpret $ \sigma(G^j_t )$  as the state-wise remaining information of $\mathcal{G}_t $ that is not already covered by  $\sigma(Z_t) $ upon $Z_t=j$ in accordance with the following result.
\begin{lemma}\label{GZG}
For $j \in S$ it holds that
\begin{align*}
\sigma( \mathcal{Z}_t,\zeta_{\mathcal{Z}_t}) \cap \{Z_t=j\} &= \sigma(G^j_t ) \cap \{Z_t=j\}, \quad t\geq0, \\
\sigma( \mathcal{Z}_{t-},\zeta_{\mathcal{Z}_{t-}}) \cap \{Z_{t-}=j\} &= \sigma(G^j_{t-} ) \cap \{Z_{t-}=j\}, \quad t>0.
\end{align*}
\end{lemma}
\begin{proof}
Since $\{Z_t=j\}\cap \{J^j_t=n\} = \{T_{(j,n)} \leq t < S_{(j,n)}\}$ and $\{Z_{t-}=j\}\cap\{J^j_{t-}=n\} =\{T_{(j,n)} < t \leq S_{(j,n)}\}$, the statement follows from~\eqref{GcondSigmaZeta}.
\end{proof}
For $G_t := \sum_{j\in S} I_j(t) \,G^j_t$, $t\geq0$, which defines another jump process on $E$, Lemma~\ref{GZG} implies the following Markovian representation of the reduced information $\mathcal{G}$:
\begin{align*}
  \mathcal{G}_t &= \sigma(Z_t,G_t) \vee \mathcal{A}, \quad t \geq 0,\\
  \mathcal{G}_{t-} &= \sigma(Z_{t-},G_{t-}) \vee \mathcal{A}, \quad t > 0.
\end{align*}

\begin{remark}\label{rmk:RepresentationAccordingToChristiansen2020}
Since
\begin{align}\label{RepresentationAccordingToChristiansen2020}\begin{split}
      \mathcal{G}_t &=  \sigma\Big( \{T_x \leq t < S_x \} \cap \{ \zeta_x \in A\}  :  x \in \mathcal{S}, A\in   \mathfrak{B}(E')\Big) \vee \mathcal{A}, \quad t \geq 0,\\
       \mathcal{G}_{t-} &= \sigma\Big( \{T_x < t \leq S_x \} \cap \{ \zeta_x \in A\}  :  x \in \mathcal{S}, A\in   \mathfrak{B}(E')\Big) \vee \mathcal{A}, \quad t > 0,
\end{split}\end{align}
our non-monotone information structures conform with those of~\citet{Christiansen2020}.\demormk
\end{remark}

Our definition of $\mathcal{G}$ implies that
\begin{align}\label{ass:Zt_is_Gt_measurable}
  \sigma(Z_t) \subseteq \sigma(\mathcal{Z}_t) \subseteq \mathcal{G}_t, \quad t\geq 0,
\end{align}
so the current states of $Z$ and $\mathcal{Z}$ are always observable, even under the non-monotone information setting.
\begin{lemma}\label{LemmaF=G}
We may set $\zeta_x =\xi_x$, $x \in \mathcal{S}$, and then
\begin{itemize}
\item $(G_t)_{t\geq0}=(F_t)_{t\geq0}$ and $(\mathcal{G}_t)_{t\geq0}=(\mathcal{F}_t)_{t\geq0}$,
\item $(G_{t-})_{t>0} = (F_{t-})_{t>0}$ and $(\mathcal{G}_{t-})_{t>0} = (\mathcal{F}_{t-})_{t>0}$.
\end{itemize}
\end{lemma}
\begin{proof}
The values of $\xi_x$ can be embedded into a separable completely metrizable space $E'$ since $\mathcal{S}$ is countable.
\end{proof}

\section{Stochastic retirement: Setup}\label{sec:setup_stoc_re}

\begin{figure}[ht]
	\centering
	\scalebox{0.75}{
	\begin{tikzpicture}[node distance=2em and 0em]
		\node[punkt] (2) {disabled};
		\node[anchor=north east, at=(2.north east)]{$2$};
		\node[punkt, left = 40mm of 2] (1) {active};
		\node[anchor=north east, at=(1.north east)]{$1$};
		\node[draw = none, fill = none, left = 18 mm of 2] (test) {};
		\node[punkt, below = 20mm of test] (3) {retired};
		\node[anchor=north east, at=(3.north east)]{p};		
		\node[punkt, below = 20mm of 3] (4) {dead};
		\node[anchor=north east, at=(4.north east)]{d};
	\path
		(1)	edge [pil, bend left = 15]				node [above]			{}		(2)
		(2)	edge [pil, bend left = 15]				node [above]			{}		(1)
		(1)	edge [pil]							node [above]			{}		(3)
		(2)	edge [pil]							node [above]			{}		(3)
		(3)	edge [pil]							node [above]			{}		(4)
		(1)	edge [pil, in = 180, out = 270]			node [above]			{}		(4)
		(2)	edge [pil, in = 0, out = 270]				node [above]			{}		(4)
	;
	\end{tikzpicture}}
	\caption{Extension of the classic three-state disability model with recovery to allow for stochastic retirement.}
	\label{fig:4state}
\end{figure}
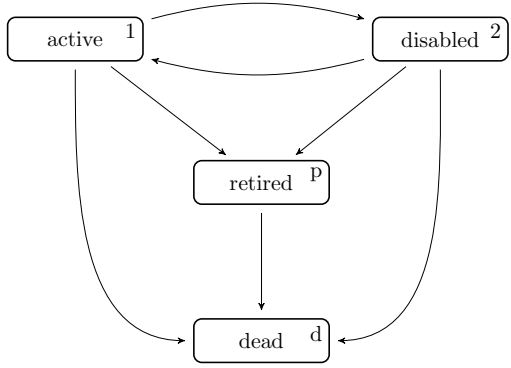

In this section, and also in Sections~\ref{sec:examples}--\ref{sec:examples2}, we consider a finite state space $S=\{1, \ldots, \varsigma, \text{p}, \text{d}\}$, where $1, \ldots, \varsigma$, $\varsigma\in\mathbb{N}$, describe different statuses of health before retirement,  $\text{p}$ means that the insured is retired and alive (pensioner, retiree), and $\text{d}$ is an absorbing state taken when the insured dies. Let $\delta$ and $\eta$ be the time of death and  the time of retirement, respectively:
\begin{align*}
\delta &:=  \inf\{ t \geq 0 : Z_t = \text{d} \}, \\
\eta &:= \inf\{ t \geq 0 : Z_t = \text{p} \}.
\end{align*}
We also here use the convention $\inf \emptyset = \infty$. Once the insured has retired, the only possible next jump shall be to state $\text{d}$, i.e.\ the counting processes $N_{\text{p} k}$, $k\in\{1,\ldots,\varsigma\}$, are almost surely constantly zero. Note that the structure of the state space entails that the insurer is not updating the health status of the insured upon or after retirement. In Figure~\ref{fig:4state} we have exemplified this setup for the case $\varsigma=2$ corresponding to a disability model allowing for recovery before retirement.

In actuarial practice, it is not uncommon to impose some restrictions on the intertemporal dependence structure of the random pattern of states $Z$ by assuming it to be e.g.\ Markovian or semi-Markovian. Under such assumptions, one may assess the mortality of retirees (pensioners) without any regard to their health status before retirement resulting in a coarser pooling of risks. If one assumes $Z$ to be a Markov process, this is related to the use of a standard mortality table such as the longevity benchmark of the Danish Financial Supervisory Authority, cf.~\citet{JarnerMoller2015}. The availability of information pertaining to the health status before retirement might depend on to what extent the insured has had health and disability coverage, and therefore it might differ between insured. In other words, the value of $\varsigma$ and the interpretation of the states $\{1,\ldots,\varsigma\}$ might not be identical among insured. If the assumption of Markovianity or semi-Markovianity is unreasonable, the above facts, and maybe even certain deliberations regarding actuarial fairness, incentivize taking into account information discarding upon and after retirement. In the following, we illustrate the potential insufficiency of (semi-)Markovian assumptions and show how to analyze the situation from the point of view of non-mononotone information. This leads to two non-trivial examples of non-monotonicity resulting from information discarding upon and after retirement.

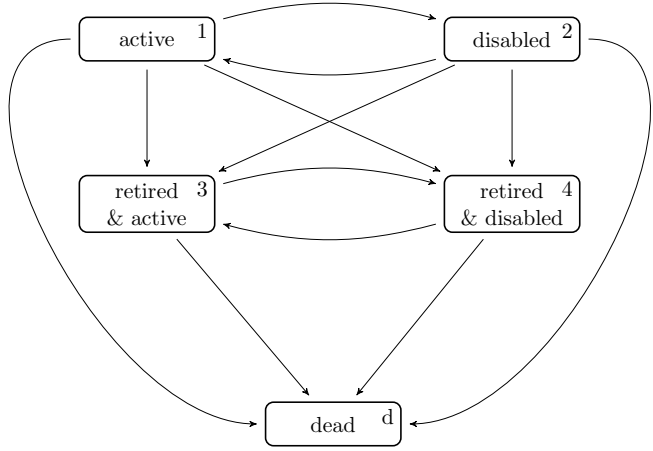
\begin{figure}[ht]
	\centering
	\scalebox{0.75}{
	\begin{tikzpicture}[node distance=2em and 0em]
		\node[punkt] (2) {disabled};
		\node[anchor=north east, at=(2.north east)]{$2$};
		\node[punkt, left = 40mm of 2] (1) {active};
		\node[anchor=north east, at=(1.north east)]{$1$};
		\node[draw = none, fill = none, left = 18 mm of 2] (test) {};
		\node[draw = none, fill = none, below = 20 mm of test] (test2) {};
		\node[punkt, below = 20mm of 2] (4) { retired \\ \& disabled};
		\node[anchor=north east, at=(4.north east)]{$4$};
		\node[punkt, below = 20mm of 1] (3) {retired \& active};
		\node[anchor=north east, at=(3.north east)]{$3$};		
		\node[punkt, below = 40mm of test2] (5) {dead};
		\node[anchor=north east, at=(5.north east)]{d};
	\path
		(1)	edge [pil, bend left = 15]				node [above]			{}		(2)
		(2)	edge [pil, bend left = 15]				node [above]			{}		(1)
		(3)	edge [pil, bend left = 15]				node [above]			{}		(4)
		(4)	edge [pil, bend left = 15]				node [above]			{}		(3)
		(1)	edge [pil]							node [above]			{}		(3)
		(2)	edge [pil]							node [above]			{}		(3)
		(1)	edge [pil]							node [above]			{}		(4)
		(2)	edge [pil]							node [above]			{}		(4)		
		(3)	edge [pil]							node [above]			{}		(5)
		(4)	edge [pil]							node [above]			{}		(5)		
		(1)	edge [pil, in = 180, out = 180]			node [above]			{}		(5)
		(2)	edge [pil, in = 0, out = 0]				node [above]			{}		(5)
	;
	\end{tikzpicture}}
	\caption{Extension of the model from Figure~\ref{fig:4state} to allow observation of the health status of the insured also upon and after retirement.}
	\label{fig:5state}
\end{figure}

It is natural to imagine the random pattern of states $Z$ as embedded into a larger fictional setting, where the insurer continues to observe the health status of the insured even upon and after retirement.  Let $\tilde{Z}$ be a random pattern of states on an extended state space $\tilde{S} := \{1,\ldots,\varsigma+1,\ldots,2\varsigma,\text{d}\}$, and let $\tilde{N}_{jk}$, $j,k\in\tilde{S}$, $j\neq k$, be the counting processes associated with $\tilde{Z}$. The states $\varsigma+1, \ldots, 2\varsigma$ replace the former state $\text{p}$ and describe different health statuses upon and after retirement. The state $\text{d}$ remains absorbing. We assume yet again that upon retirement, the insured stays retired or dies, i.e.\ the counting processes $\tilde{N}_{jk}$, $j\in\{\varsigma + 1, \ldots, 2\varsigma\}$, $k\in\{1,\ldots,\varsigma\}$, are almost surely constantly zero. In Figure~\ref{fig:5state} we have exemplified this setup for the case $\varsigma=2$ corresponding to a disability model allowing for recovery before and after retirement.

In similar fashion as~\eqref{assump:tech_noexplo}, we suppose that
\begin{align}\label{assump:tech_noexplo_tilde}
  \E\bigg[\sum_{j,k\in \tilde{S}, j\neq k}\! \tilde{N}_{jk}(t) \,\bigg] < \infty, \quad t\geq0.
\end{align}
This in particular ensures that $\tilde{Z}$ is non-explosive, and we may recover $Z$ from $\tilde{Z}$ according to
\begin{align*}
Z_t =
\begin{cases}
\tilde{Z}_t & \text{ if }\tilde{Z}_t \in \{1,\ldots,\varsigma,\text{d}\}, \\
\text{p} & \text{ if }\tilde{Z}_t \in \{\varsigma+1,\ldots,2\varsigma\}
\end{cases}
\end{align*}
for all $t \geq 0$.

The information available to the insured is represented by the completed natural filtration $\tilde{\mathcal{F}}=(\tilde{\mathcal{F}}_t)_{t\geq}$ of $\tilde{Z}$ with components given by
\begin{align*}
\tilde{\mathcal{F}}_t = \sigma( \tilde{Z}_s : s\leq t) \vee \mathcal{A}
\end{align*}
for $t\geq0$. In most instances, the insurer must resort to the information $\mathcal{F}$ rather than the extended information $\tilde{\mathcal{F}}$, since disability and health coverage typically ceases upon retirement.

It appears consistent with actuarial practice to propose that the underlying random pattern of states $\tilde{Z}$ is Markovian or semi-Markovian. We now study the resulting implications on $Z$. Let $U=(U_t)_{t\geq 0}$ and $\tilde{U}=(\tilde{U}_t)_{t \geq 0}$ be the duration processes associated with $Z$ and $\tilde{Z}$, respectively, given by
\begin{align*}
U_t & := t - \sup\{s \in [0,t] : Z_s \neq Z_t\}, \\
\tilde{U}_t & := t - \sup\{s \in [0,t] : \tilde{Z}_s \neq \tilde{Z}_t\},
\end{align*}
for $t\geq0$. Note that $\mathbf{1}_{\{t\leq\eta\}}U_t = \mathbf{1}_{\{t\leq\eta\}}\tilde{U}_t$. Let $U^{\text{p}}=(U^{\text{p}}_t)_{t \geq 0}$ be the time since retirement given by
\begin{align*}
U^{\text{p}}_t
:=
\begin{cases}
0 & \text{ if } 0 \leq t < \eta, \\
t - \eta & \text{ if } t \geq \eta,
\end{cases}
\end{align*}
let $H=(H_t)_{t \geq 0}$ be the state of the insured just before retirement given by
\begin{align*}
H_t :=
\begin{cases}
Z_t & \text{ if } 0 \leq t < \eta, \\
Z_{\eta-} & \text{ if } t \geq \eta,
\end{cases}
\end{align*}
with $Z_{0-} := Z_0$, and let $U^{\text{h}}=(U^{\text{h}}_t)_{t \geq 0}$ be the duration of the latest sojourn before retirement given by
\begin{align*}
U^{\text{h}}_t :=
\begin{cases}
U_t & \text{ if } 0 \leq t < \eta, \\
U_{\eta-} & \text{ if } t \geq \eta,
\end{cases}
\end{align*}
with $U_{0-} := 0$. In the following, $U_{0-}^{\text{h}} := 0$, $H_{0-} := 0$, and $U_{0-}^{\text{p}} := 0$.
\begin{proposition}\label{prop:linkmarkov}Assume that $\tilde{Z}$ satisfies~\eqref{assump:tech_noexplo_tilde}.
\begin{enumerate}
  \item[(a)] If $(\tilde{Z},\tilde{U})$ is Markovian, then $(Z,U^{\emph{h}},U^{\emph{p}},H)$ is Markovian.
  \item[(b)] If $\tilde{Z}$ is Markovian, then $(Z,U^{\emph{p}},H)$ is Markovian.
\end{enumerate}
\end{proposition}
\begin{proof}
The only non-trivial statement of the proposition relates to the intertemporal dependence structure of $Z$ after retirement, so it suffices to study the quantities
\begin{align*}
P\!\left(Z_s = \text{p} \left| \, \mathcal{F}_t \right.\right)
\end{align*}
on the event $\{Z_t = \text{p}\}$ for $0 \leq t < s < \infty$. To this end, fix such $t,s$ and consider sets $A^n_t := \{Z_t = \text{p}, \dot{N}(t) = n\}$, $n \in \mathbb{N}_0$, where $\dot{N}=(\dot{N}(t))_{t\geq0}$ is the process counting the total number of jumps of $Z$ given by $\dot{N} := \sum_{j,k\in S, j\neq k} N_{jk}$. Furthermore, denote with $\tau=(\tau_i)_{i\in\mathbb{N}_0}$ and $\tilde{\tau}=(\tilde{\tau}_i)_{i\in\mathbb{N}_0}$ the point processes corresponding to the jump times of $Z$ and $\tilde{Z}$, respectively; here $\tau_0 := 0$ and $\tilde{\tau}_0 := 0$.

On $A^n_t$ it holds that
\begin{equation}\label{eq:help_result_proof_ex}
  \begin{alignedat}{2}
&\tau_n = \tilde{\tau}_n = \eta, \quad \quad
&&\tilde{Z}_{\tilde{\tau}_n} \in \{\varsigma+1,\ldots,2\varsigma\}, \\
&\tau_i = \tilde{\tau}_i,
&&Z_{\tau_i} = \tilde{Z}_{\tilde{\tau}_i} \in \{1,\ldots,\varsigma\}, \quad i=0,\ldots,n-1.
  \end{alignedat}
\end{equation}
In particular,
\begin{align*}
P\!\left(Z_s = \text{p} \left| \, \mathcal{F}_t \right.\right)\!\mathbf{1}_{\!A^n_t}
=
P\!\left(\tilde{Z}_s \in \{\varsigma+1,\ldots,2\varsigma\} \left| \, \tilde{\tau}_n, Z_{\tilde{\tau}_n},\tilde{\tau}_{n-1},\tilde{Z}_{\tilde{\tau}_{n-1}},\ldots,\tilde{\tau}_1,\tilde{Z}_{\tilde{\tau}_1},\tilde{\tau}_0,\tilde{Z}_{\tilde{\tau}_0} \right.\right)\!\mathbf{1}_{\!A^n_t}
\end{align*}
almost surely, cf.\ Lemma~\ref{lemma:hotfix_completion_condexp} below.

We now focus on case (a), so suppose $(\tilde{Z},\tilde{U})$ is Markovian such that $\tilde{Z}$ is semi-Markovian. By the law of iterated expectations and the strong Markov property, cf.~Theorem 7.5.1 in~\citet{Jacobsen2006}, it follows that
\begin{align*}
&P\!\left(\tilde{Z}_s \in \{\varsigma+1,\ldots,2\varsigma\} \left| \, \tilde{\tau}_n, Z_{\tilde{\tau}_n},\tilde{\tau}_{n-1},\tilde{Z}_{\tilde{\tau}_{n-1}},\ldots,\tilde{\tau}_1,\tilde{Z}_{\tilde{\tau}_1},\tilde{\tau}_0,\tilde{Z}_{\tilde{\tau}_0} \right.\right)\!\mathbf{1}_{\!A^n_t} \\
&=
\,\E\!\left[\left. P\!\left(\tilde{Z}_s \in \{\varsigma+1,\ldots,2\varsigma\} \left| \, \tilde{\tau}_n, \tilde{Z}_{\tilde{\tau}_n},\tilde{\tau}_n - \tilde{\tau}_{n-1} \right.\right) \, \right| \tilde{\tau}_n, Z_{\tilde{\tau}_n},\tilde{\tau}_{n-1},\tilde{Z}_{\tilde{\tau}_{n-1}},\ldots,\tilde{\tau}_1,\tilde{Z}_{\tilde{\tau}_1},\tilde{\tau}_0,\tilde{Z}_{\tilde{\tau}_0}\right]\! \mathbf{1}_{\!A^j_t}  \\
&=
\E\!\left[\left. P\!\left(\tilde{Z}_s \in \{\varsigma+1,\ldots,2\varsigma\} \left| \, \tilde{\tau}_n, \tilde{Z}_{\tilde{\tau}_n},\tilde{\tau}_n - \tilde{\tau}_{n-1} \right.\right) \, \right| \tilde{\tau}_n, Z_{\tilde{\tau}_n},\tilde{\tau}_{n-1},\tilde{Z}_{\tilde{\tau}_{n-1}}\right]\! \mathbf{1}_{\!A^j_t} \\
&=
\,P\!\left(\tilde{Z}_s \in \{\varsigma+1,\ldots,2\varsigma\} \left| \, \tilde{\tau}_n, Z_{\tilde{\tau}_n},\tilde{\tau}_{n-1},\tilde{Z}_{\tilde{\tau}_{n-1}}\right.\right)\!\mathbf{1}_{\!A^n_t}
\end{align*}
almost surely. Thus on $A^n_t = \{Z_t = \text{p}, \dot{N}(t) = n\}$ it almost surely holds that
\begin{align*}
P\!\left(Z_s = \text{p} \left| \, \mathcal{F}_t \right.\right)
&=
P\!\left(Z_s = \text{p} \left| \, \tilde{\tau}_n, Z_{\tilde{\tau}_n},\tilde{\tau}_{n-1},\tilde{Z}_{\tilde{\tau}_{n-1}}\right.\right) \\
&=
P\!\left(Z_s = \text{p} \left| \, t - \tilde{\tau}_n, Z_{\tilde{\tau}_n},t - \tilde{\tau}_{n-1},\tilde{Z}_{\tilde{\tau}_{n-1}}\right.\right) \\
&=
P\!\left(Z_s = \text{p} \left| \,  U^{\text{p}}_t, Z_t,U^{\text{h}}_t,H_t\right.\right)\!,
\end{align*}
cf.~\eqref{eq:help_result_proof_ex}. Since the final line does not depend on $n\in\mathbb{N}_0$, and since $\{Z_t = \text{p}, \dot{N}(t) = \infty\}$ is a null-set according to Lemma~\ref{FinitelyManyJumps}, we conclude that if $(\tilde{Z},\tilde{U})$ is Markovian, then on $\{Z_t = \text{p}\}$ we almost surely have
\begin{align*}
P\!\left(Z_s = \text{p} \left| \, \mathcal{F}_t \right.\right)
=
P\!\left(Z_s = \text{p} \left| \,  U^{\text{p}}_t, Z_t,U^{\text{h}}_t,H_t\right.\right)
\end{align*}
proving (a). The proof of (b) follows by similar arguments.
\end{proof}

It is possible to derive necessary and sufficient conditions for which (semi-)Markovianity of $\tilde{Z}$ implies (semi-)Markovianity of $Z$, see e.g.~\citet{Serfozo1971}. In general, such conditions are very restrictive and do not apply to models of actuarial relevance. So we can normally expect that if $(\tilde{Z},\tilde{U})$ is Markovian, then due to Proposition~\ref{prop:linkmarkov}(a), $(Z,U)$ is not Markovian, and that if $\tilde{Z}$ is Markovian, then due to Proposition~\ref{prop:linkmarkov}(b), $Z$ is not Markovian and most likely $(Z,U)$ is also not Markovian. Consequently, the possible lack of availability of information pertaining e.g.\ to the health status before retirement is not trivially resolved but has to be taken explicitly into account. This can be done by analyzing the situation from the point of view of non-monotone information. To this end, we introduce two sub-sequences of $\sigma$-algebras $\mathcal{G}^{(1)}=(\mathcal{G}^{(1)}_t)_{t\geq 0}$ and $\mathcal{G}^{(2)}=(\mathcal{G}^{(2)}_t)_{t\geq 0}$ given by
\begin{align*}
\mathcal{G}^{(1)}_t
&:=
\sigma(Z_s \mathbf{1}_{\{Z_t \in \{1,\ldots,\varsigma\}\}}, \mathbf{1}_{\{\eta \leq s\}}, \mathbf{1}_{\{\delta \leq s\}} : s \leq t) \vee \mathcal{A}, \\
\mathcal{G}^{(2)}_t
&:=
\sigma(Z_s \mathbf{1}_{\{Z_t \in \{1,\ldots,\varsigma\}\}}, \mathbf{1}_{\{\eta \leq t\}}, \mathbf{1}_{\{\delta \leq s\}} : s \leq t) \vee \mathcal{A}
\end{align*}
for $t\geq0$. The information $\mathcal{G}^{(1)}$ corresponds to the case where upon retirement or death the insurer does not have access to or discards all previous records regarding the health status of the insured. The sub-sequence $\mathcal{G}^{(2)}$ of $\mathcal{G}^{(1)}$ even keeps no record of the time of retirement. Note that we also discard past information upon death, but for most if not all practical purposes this does not matter.

In Sections~\ref{sec:examples}--\ref{sec:examples2}, we study the dynamics of prospective reserves under information $\mathcal{G}^{(1)}$ and $\mathcal{G}^{(2)}$, when the following link to the general setting becomes key.

\begin{lemma}\label{lemma:translation}
The $\sigma$-algebras $\mathcal{G}_t^{(1)}$  and $\mathcal{G}_{t}^{(2)}$,
$t \geq 0$,  can be brought on the form of Section~\ref{subsec:nonmono}.
\end{lemma}
\begin{proof}
We obtain $\mathcal{G}_t^{(1)} =\sigma(Z_t,G_t) \vee \mathcal{A}$
for
\begin{alignat*}{3}
 \zeta_{(\text{p},0)} &:= \eta ,\quad  \zeta_{(\text{d},0)} &:=(\eta, \delta), \quad \zeta_{(j,n)} &:= \xi_{(j,n)}, \quad j \in \{1, \ldots, \varsigma\}, \, n \in \mathbb{N}_0,
\end{alignat*}
and we obtain $\mathcal{G}_t^{(2)} =\sigma(Z_t,G_t) \vee \mathcal{A}$
for
\begin{alignat*}{3}
 \zeta_{(\text{p},0)} &:= 0 ,\quad
 \zeta_{(\text{d},0)} &:=\delta, \quad
 \zeta_{(j,n)} &:= \xi_{(j,n)}, \quad j \in \{1, \ldots, \varsigma\}, \, n \in \mathbb{N}_0.
\end{alignat*}
Note that since $Z$ almost surely visits each of the states $\text{p}$ and $\text{d}$ at most once, the definition of $\zeta_{\text{p},n}$ as well as $\zeta_{\text{d},n}$ for $n\in\mathbb{N}$ is inconsequential.
\end{proof}
From this point onward, for $\mathcal{G}^{(1)}$ and $\mathcal{G}^{(2)}$ the random variables $ \zeta_x$, $x \in  \mathcal{S}$, are always taken to be those from the proof of Lemma~\ref{lemma:translation}.

\section{State-wise prospective reserves}\label{sec:state_wise_prores}

We now turn our attention to the concept of state-wise prospective reserves. A life insurance contract between the insured and the insurer is stipulated by the specification of a payment process $B=(B(t))_{t \geq 0}$ representing the accumulated benefits minus premiums.
 We generally suppose that
\begin{align}\label{eq:nice_payments}
B(\d t)&=  \sum_{j\in S} I_j(t-) \, b_j(t) \, \beta (\d t)
+
\sum_{j,k\in S\atop j \neq k} b_{jk}(t) \, N_{jk}(\d t),
\end{align}
where $b_j$ and $b_{jk}$ are $\mathcal{F}$-predictable bounded processes and the measure $\beta$ is a sum of the Lebesgue-measure $m$ and a countable number of Dirac-measures $(\varepsilon_{t_n})_{n\in\mathbb{N}}$:
\begin{align*}
  \beta (A):=  m(A) + \sum_{n=1}^{\infty} \varepsilon_{t_n}(A), \quad A \in \mathfrak{B}([0,\infty)),
\end{align*}
for deterministic time points $0 \leq t_1 < t_2 < ... $ that are increasing to infinity (i.e.~there are at most a finite number of such time points on each compact interval). By using the measure $\beta$ our model setting can simultaneously cover absolutely continuous as well as most types of discrete sojourn payments. In the following, we deal with finite variation processes, i.e.\ adapted processes with c\`{a}dl\`{a}g paths of finite variation on compacts. For this class of processes, stochastic integrals equal path-wise Lebesgue integrals. To both stress this point of view and to avoid confusion when several arguments are involved, we have replaced and are going to continue to replace the notation $\d X(t)$ by $X(\mathrm{d}t)$.

Consider a deterministic bank account $\kappa : [0,\infty) \mapsto (0,\infty)$ assumed measurable, bounded away from zero,  c\`{a}dl\`{a}g, and of finite variation on compacts, with initial value $\kappa(0)=1$. Denote with $v$ the corresponding discount function given by
\begin{align*}
[0,\infty) \ni t \mapsto v(t) = \frac{1}{\kappa(t)}.
\end{align*}
Let $Y=(Y(t))_{t\geq0}$ give the discounted accumulated future payments, formally defined as
\begin{align*}
Y(t):=
\int_{(t,\infty)} \frac{\kappa(t)}{\kappa(s)}\,  B(\d s), \quad t\geq0.
\end{align*}
We suppose that $Y$ has finite expected variation on compacts, which is e.g.\ the case if there exists $n\geq0$ such that $B(t) = B(n)$ for $t > n$. Note that $Y$ then in particular is a finite variation process. The prospective reserve (under information $\mathcal{G}$) is denoted $V=(V(t))_{t\geq0}$ and defined as the optional projection of $Y$ (w.r.t.\ $\mathcal{G}$) satisfying
\begin{align*}
V(t) = \E[ Y(t) \,  | \,  \mathcal{G}_{t}]
\end{align*}
almost surely for each $t\geq0$. The prospective reserve is unique (up to evanescence), and since $Y$ has finite expected variation on compacts, the prospective reserve $V$ is actually a finite variation process, cf.\ Theorem 4.1 in~\citet{Christiansen2020}.

The following result shows that whenever we consider conditional expectations, the null-sets $\mathcal{A}$ can be disregarded. This fact is used throughout the paper to allow us to work with the information $\sigma(Z_t,G_t)$ instead of the information $\mathcal{G}_t$.
\begin{lemma}\label{lemma:hotfix_completion_condexp}
Let $X$ be an integrable random variable and let $\mathcal{H}$ be some
sub-$\sigma$-algebra of $\mathbb{F}$. Then
$\E[X \,| \, \mathcal{H}]=\E[X \, | \, \mathcal{H}\vee \mathcal{A}]$ almost surely.
\end{lemma}
\begin{proof}
The random variable $\E[X \, | \, \mathcal{H}]$ satisfies the defining
Radon-Nikodym equation of the conditional expectation $\E[X \, | \,
\mathcal{H}\vee \mathcal{A}]$.
\end{proof}

\subsubsection*{General definition and explicit construction}

From the perspective of the Thiele differential equations and the stochastic Thiele equation,  state-wise prospective reserves should be processes  $(V_j)_{j\in S}$ such that
\begin{itemize}
  \item $  V(t) = \sum_{j \in S} I_j(t)\, V_j(t)$,
  \item  $V_j(t)$ does not depend on the actual value of $Z_t$,
  \item the Ito differentials $I_j(t-)\, \d V_j(t)$ exist.
\end{itemize}
This motivates the following definition:

\begin{definition}\label{def:stat_counters}
We denote  $(V_j(t))_{t\geq 0,j\in S}$ as \textit{state-wise prospective reserves}  if for each $j \in S$
\begin{enumerate}
  \item[(1)] $I_j(t) \, V_j(t) = I_j(t)\, V(t)$ almost surely for each $t \geq 0$,
  \item[(2)] $I_j(t-) \, V_j(t)$ is $\mathcal{G}_{t-}$-measurable for each $t >0$,
  \item[(3)] the mapping $t \mapsto \int_{(0,t]} I_j(s-)\, V_j(\d s)$ defines a finite variation process.\demodef
\end{enumerate}
\end{definition}
Because of Lemma~\ref{GZG}, properties (2) and (1) yield
\begin{align*}
I_j(t-)\, V_j(t)  \E[I_j(t) \,| \,  G^j_{t-}, Z_{t-}=j] &= I_j(t-) \E[ I_j(t) V_j(t)  \,|\,  \mathcal{G}_{t-}]  \\
&= I_j(t-) \E[ I_j(t) V(t)  \,|\,  \mathcal{G}_{t-}]\\
&= I_j(t-) \E[ I_j(t) V(t)  \,|\,  G^j_{t-}, Z_{t-}=j]
\end{align*}
almost surely. Consequently,
\begin{align*}
I_j(t-)\, V_j(t) &=  I_j(t-) \frac{\E[ I_j(t) V(t)  \, | \,  G^j_{t-}, Z_{t-}=j]}{ \E[I_j(t) \,| \,  G^j_{t-}, Z_{t-}=j]}
\end{align*}
almost surely unless the denominator is zero. The latter equation suggests to set $V_j(t) = \E[ I_j(t) V(t) \, | \,  G^j_{t-}, Z_{t-}=j]/ \E[I_j(t) \,| \,  G^j_{t-}, Z_{t-}=j]$. As this quantity is
$\sigma(G^j_{t-})$-measurable, we can alternatively derive the point-wise almost sure equation $V_j(t) \,  \E[I_j(t) \, | \, G^j_{t-}] =\E[I_j(t)\, V(t) \, | \, G^j_{t-}]$, which almost surely implies
\begin{align}\label{IntuitDef}
V_j(t)=
\frac{
\E[I_j(t)\, V(t)  \,|\, G^j_{t}]}{\E[I_j(t) \,|\, G^j_{t}]}
\end{align}
because of $G^j_{t-}=G^j_t$, unless the denominator is zero. The denominator is non-zero on $\{Z_t=j\}$, since $\E[I_j(t) \mathbf{1}_{\{\E[I_j(t) \, | \,  G^j_{t}]=0\}}]=\E[\E[I_j(t) \, | \,  G^j_{t}] \mathbf{1}_{\{\E[I_j(t) \,|\, G^j_{t}]=0\}}]=0$. So \eqref{IntuitDef} is at least well-defined on $\{Z_t=j\}$, which is sufficient for property (1) in Definition \ref{def:stat_counters}.
 But the construction \eqref{IntuitDef} offers a  point-wise almost sure definition only, which means that $V_j$ is not necessarily well-defined as a stochastic process and no path properties can be derived. In order to solve that issue, we temporarily switch from $Z$ to the expanded state process $\mathcal{Z}$ and exploit the fact that
\begin{align*}
\sigma(G^j_t) \cap   \{\mathcal{Z}_t=(j,n)\} =  \sigma( \zeta_{(j,n)} )  \cap  \{\mathcal{Z}_t=(j,n)\}, \quad j\in S, n\in\mathbb{N}_0,
\end{align*}
see Lemma \ref{GZG} and \eqref{GcondSigmaZeta}. The advantage of the right hand side is that $\sigma( \zeta_{(j,n)} )$ is constant in time, whereas $\sigma(G^j_t)$ changes with $t$.

Under the convention $0/0:=0$, let
\begin{align}\label{DefOfVx}
\mathcal{V}_x(t) :=\frac{\E[ \mathcal{I}_x(t) Y(t) \, | \, \zeta_x]}{\E[ \mathcal{I}_x(t) \,|\, \zeta_x]}, \quad x\in\mathcal{S}, t\geq0.
\end{align}
This definition is similar to~\eqref{IntuitDef}, but here the conditions are constant in time, so $\mathcal{V}_x$ is well-defined as a stochastic process.
\begin{proposition}\label{Proposition:RegVx}
The mapping $t \mapsto \int_{(0,t]} \mathcal{I}_x(s-) \,  \mathcal{V}_x(\d s)$ defines a finite variation process for each $x \in \mathcal{S}$.
  For each $t \geq 0 $ and $x \in \mathcal{S}$ we almost surely have
\begin{align}\label{PRE:DefProspResGExt}
  \mathcal{V}_x(t) = \E[ Y(t)  \,|\, \zeta_{x}, \mathcal{Z}_t=x ] \quad \text{ on } \{\mathcal{Z}_t=x\}.
\end{align}
\end{proposition}
Recall that the conditional expectation in \eqref{PRE:DefProspResGExt} is a short-hand notation for the
factorized conditional expectation $\E[ Y(t) \, | \,  (\zeta_{x}, \mathcal{Z}_t)=\cdot \,]$  evaluated at $(\zeta_x, x)$. This kind of short-hand notation is used throughout the paper.

\begin{proof}
According to Remark~\ref{rmk:RepresentationAccordingToChristiansen2020}, we can use the results from~\citet{Christiansen2020} in order to verify the statements of the proposition. By following the argumentation in the proof of Theorem 4.1 in \citet{Christiansen2020}, we can show that $t \mapsto \mathbf{1}_{\{\mathcal{Z}_{t}=x\}}\, \mathcal{V}_x(t)$ has c\`{a}dl\`{a}g paths of finite variation on compacts.
The latter path properties do not only hold on  $\{\mathcal{Z}_t=x \} =\{T_x \leq t < S_x\} $ but even on $\{\mathcal{Z}_t=x \textrm{ or }  \mathcal{Z}_{t-}=x \}=\{T_x \leq t \leq S_x\}$ because of Lemma \ref{FinitelyManyJumps}. The process $t \mapsto \int_{(0,t]} \mathcal{I}_x(s-) \,  \mathcal{V}_x(\d s)$ is indeed $\mathcal{F}$-adapted because of \eqref{GcondSigmaZeta}. Proposition 4.2 in~\citet{Christiansen2020} and~\eqref{DefG} yield~\eqref{PRE:DefProspResGExt}.
\end{proof}
Proposition \ref{Proposition:RegVx} implies that
\begin{align*}
   \mathcal{V}(t)=\sum_{x \in \mathcal{S}} \mathcal{I}_x(t)\,  \mathcal{V}_x(t)
\end{align*}
almost surely for each $t\geq 0$. Thus,  we can think of $\mathcal{V}_x(t)$, $x\in\mathcal{S}$, as state-wise prospective reserves with respect to the expanded state process $\mathcal{Z}$. The following theorem shows how to obtain state-wise prospective reserves according to Definition~\ref{def:stat_counters} starting from such quantities.
\begin{theorem}\label{Proposition:RegVj}
The processes $(V_j)_{j\in S}=(V_j(t))_{t\geq0,j\in S}$  defined by
\begin{align}\label{RegDefOfVj}
  V_j(t) := \sum_{n=0}^{\infty} \mathbf{1}_{\{J^j_t=n\}} \mathcal{V}_{(j,n)}(t), \quad t\geq0,
\end{align}
 are state-wise prospective reserves in the sense of Definition~\ref{def:stat_counters}. If   $(\widetilde{V}_j)_{j\in S}$ is another set of state-wise prospective reserves, then $(I_j\widetilde{V}_j)_{j\in S}$ almost surely equals $(I_jV_j)_{j\in S}$.
\end{theorem}
\begin{proof}
By applying Proposition \ref{Proposition:RegVx} and using that $(\mathcal{Z}_t,\zeta_{\mathcal{Z}_t})= (Z_t,G_t)$, we get that
\begin{align}\label{DefProspResG}
 I_j(t)\, V_j(t) =   I_j(t) \sum_{n=0}^{\infty} \mathbf{1}_{\{J^j_t=n\}}  \E[ Y(t) \, | \, \zeta_{(j,n)}, \mathcal{Z}_t=(j,n) ]  =  I_j(t)\,\E[ Y(t) \,  | \, Z_t=j,G^j_t]
  \end{align}
almost surely for each $t \geq 0 $ and $j \in S$. This fact and Lemma \ref{GZG} imply $I_j(t)\, V_j(t) = I_j(t)\, V(t)$ almost surely.  Moreover, the $\mathcal{G}_{t-}$-measurability of $I_j(t-)V_j(t)$ follows from~\eqref{DefG} and the fact that $J_{t-}^j = J_{t}^j$.
By applying Proposition \ref{Proposition:RegVx} we get that  $t \mapsto \int_{(0,t]} I_j(s-)\, V_j(\d s)$ is a finite variation process.
If   $(\widetilde{V}_j)_{j\in S}$ is another set of state-wise prospective reserves, then $I_j(t) V_j(t)=I_j(t)\widetilde{V}_j(t)$ almost surely for each $t$. By similar arguments as in the proof of Proposition~\ref{Proposition:RegVx}, property (3) of Definition~\ref{def:stat_counters} implies that $I_j\widetilde{V}_j$ and $I_jV_j$ have c\`{a}dl\`{a}g paths, so they are even almost surely equal as stochastic processes.
\end{proof}

\subsubsection*{Further state-dependent quantities}

The equations \eqref{PRE:Thiele} and \eqref{PRE_generalizedThiele} involve further kinds of state-dependent prospective reserves, which, so far, have only been defined point-wise almost surely. We now clarify their definitions as processes and introduce some additional state-dependent quantities which are needed later on.
\begin{definition}\label{DefFurtherStateWiseQuant2}
For $x,y \in \mathcal{S}$, $x\neq y$, $e \in E'$, and $t>0$ let
\begin{align*}
  \mathcal{V}_{xx}(t) &:=\frac{\E[ \mathcal{I}_x(t-)\mathcal{I}_x(t) Y(t) \, | \, \zeta_x]}{\E[ \mathcal{I}_x(t-)\mathcal{I}_x(t) \, | \, \zeta_x]},\\
\mathcal{V}_{xy}(t) &:=\frac{\E[ \mathcal{I}_x(T_y-) Y(T_y) \, | \, \zeta_x,T_y=t]}{\E[ \mathcal{I}_x(T_y-) \, | \, \zeta_x,T_y=t]},\\
  \mathcal{V}_{xy}(t,e) &:=\frac{\E[ \mathcal{I}_x(T_y-) Y(T_y) \, | \, \zeta_x,T_y=t,\zeta_y=e]}{\E[ \mathcal{I}_x(T_y-) \, | \, \zeta_x,T_y=t,\zeta_y=e]},\\
  \overline{\mathcal{V}}_{yx}(t,e) &:=\frac{\E[ \mathcal{I}_x(S_y) Y(S_y) \, | \, S_y=t,\zeta_y=e,\zeta_x]}{\E[ \mathcal{I}_x(S_y) \, | \,  S_y=t,\zeta_y=e,\zeta_x]}
\end{align*}
under the convention $0/0:=0$.\demodef
\end{definition}
Since the random quantities $\zeta_x, T_y, \zeta_y,S_y$ do not depend on $t$ and the time argument enters  the conditions  in
Definition~\ref{DefFurtherStateWiseQuant2} only through factorization, there are at most  a countable number of exception zero sets involved in the definitions. Consequently, the quantities form well-defined processes that are jointly measurable as mappings of $(\omega,t)$ and of $(\omega,t,e)$, respectively.
\begin{proposition}\label{Prop:MeaningOfVxxVxyVyx}
For each $x,y \in \mathcal{S}$, $x \neq y$, $e \in E'$, and $t>0$ we almost surely have
\begin{align*}
  \mathcal{V}_{xx}(t)&= \E[ Y(t) \, | \, \mathcal{Z}_{t-}=x , \zeta_x, \mathcal{Z}_{t}=x ],\\
  \mathcal{V}_{xy}(t)&= \E[ Y(t) \, | \, \mathcal{Z}_{t-}=x,\zeta_{x},\mathcal{Z}_{t}=y ],\\
  \mathcal{V}_{xy}(t,e)&= \E[ Y(t) \, | \, \mathcal{Z}_{t-}=x,\zeta_{x},\mathcal{Z}_{t}=y,\zeta_{y}=e ]
\end{align*}
on $\{\mathcal{Z}_{t-}=x\}$ and
 \begin{align*}
  \overline{\mathcal{V}}_{yx}(t,e)&= \E[ Y(t) \, | \, \mathcal{Z}_{t-}=y,\zeta_{y}=e,\mathcal{Z}_{t}=x,\zeta_{x} ]
\end{align*}
on $\{\mathcal{Z}_{t}=x\}$.
\end{proposition}
\begin{proof}
Apply  Proposition 4.2 in  \citet{Christiansen2020} and use that
\begin{align*}
  \{ \mathcal{Z}_{t-}=x,  \mathcal{Z}_{t}=y,\zeta_y=e\}&= \{ \mathcal{Z}_{t-}=x, T_y=t,\zeta_y=e\},\\
  \{ \mathcal{Z}_{t-}=x, \zeta_x=e, \mathcal{Z}_{t}=y\}&= \{  S_x=t,\zeta_x=e,\mathcal{Z}_{t}=y\}.\qedhere
\end{align*}
\end{proof}

Based on  Proposition \ref{Prop:MeaningOfVxxVxyVyx}, we
 can interpret the state-wise quantities $\mathcal{V}_{xx}(t)$, $\mathcal{V}_{xy}(t,e)$, $\overline{\mathcal{V}}_{yx}(t,e)$ as follows:
\begin{itemize}
\item $\mathcal{V}_{xx}(t)$ is the prospective reserve for staying in state $x$ at time $t$: If in state $x$ at time $t-$ or time $t$, what one would set aside in case no change in $\mathcal{Z}$ occurs at time $t$.
\item $\overline{\mathcal{V}}_{yx}(t,e)$ is the \textit{backward} prospective reserve at transition of $\mathcal{Z}$ from $y$ to $x$ with original information $e$ at time $t$: If in state $x$ at time $t$, what one would set aside in case a change \underline{from} state $y$ occurred and $\zeta_y=e$.
\item $\mathcal{V}_{xy}(t,e)$ is the \textit{forward} prospective reserve at transition of $\mathcal{Z}$  from $x$ to $y$ with updated information $e$ at time $t$: If in  state $x$ at time $t-$, what one would set aside in case a change \underline{to} state $y$ occurred and $\zeta_y=e$.
\end{itemize}
\begin{definition}\label{DefFurtherStateWiseQuant1}
For $j,k\in S$, $j \neq k$, $(n,e)\in E$, and $t>0$ let
\begin{align*}
  V_{jj}(t) &:= \sum_{l=0}^{\infty} \mathbf{1}_{\{J^j_{t-}=l\}} \mathcal{V}_{(j,l)(j,l)}(t),\\
   V_{jk}(t) &:= \sum_{l,n=0}^{\infty} \mathbf{1}_{\{J^j_{t-}=l\}} \mathbf{1}_{\{J^k_{t}=n\}}\mathcal{V}_{(j,l)(k,n)}(t),\\
  V_{jk}(t,(n,e)) &:= \sum_{l=0}^{\infty} \mathbf{1}_{\{J^j_{t-}=l\}} \mathcal{V}_{(j,l)(k,n)}(t,e),\\
\hspace{40mm}  \overline{V}_{kj}(t,(n,e)) &:= \sum_{l=0}^{\infty} \mathbf{1}_{\{J^j_{t}=l\}} \overline{\mathcal{V}}_{(k,n)(j,l)}(t,e).\hspace{50mm}\triangle
\end{align*}
\end{definition}

Definition~\ref{DefFurtherStateWiseQuant1} involves a countable number of processes from Definition~\ref{Prop:MeaningOfVxxVxyVyx} only, so the quantities are well-defined as stochastic processes and are jointly measurable as mappings of $(\omega,t)$ and of $(\omega,t,(n,e))$, respectively.
\begin{proposition}\label{Prop:MeaningOfVjjVjkVkj}
For each $j,k\in S$, $j\neq k$, $g \in E$, and $t>0$ we almost surely have
\begin{align*}
  V_{jj}(t)&=\E[ Y(t) \, | \, Z_{t-}=j,G^j_{t-},Z_t=j ],\\
  V_{jk}(t)&=\E[ Y(t) \, | \, Z_{t-}=j,G^j_{t-}, Z_t=k ],\\
  V_{jk}(t,g)&=\E[ Y(t) \, | \, Z_{t-}=j,G^j_{t-}, Z_t=k,G^k_t=g ]
\end{align*}
on $\{Z_{t-}=j\}$ and
\begin{align*}
  \overline{V}_{kj}(t,g)&=\E[ Y(t) \, | \,  (Z_{t-}=k,G^k_{t-}=g,Z_t=j,G^j_t ]
\end{align*}
on $\{Z_t = j\}$.
\end{proposition}
\begin{proof}
  Apply Proposition \ref{Prop:MeaningOfVxxVxyVyx} and  use that $\mathcal{Z}_t= (Z_t,J_t)$, $\mathcal{Z}_{t-}= (Z_{t-},J_{t-})$ and $(J_t,\zeta_{\mathcal{Z}_t})=G_t $, $(J_{t-},\zeta_{\mathcal{Z}_{t-}})=G_{t-} $.
\end{proof}
Based on  Proposition \ref{Prop:MeaningOfVjjVjkVkj},  we
interpret the state-wise quantities $V_{jj}(t)$, $V_{jk}(t,g)$, $\overline{V}_{kj}(t,g)$ as follows:
\begin{itemize}
\item $V_{jj}(t)$ is the prospective reserve for staying in state $j$ at time $t$.
\item $\overline{V}_{kj}(t,g)$ is the \textit{backward} prospective reserve at transition of $Z$ from $k$ to $j$ with $G^k_{t-}=g$ at time $t$.
\item $V_{jk}(t,g)$ is the \textit{forward} prospective reserve at transition of $Z$ from $j$ to $k$ with $G^k_t=g$ at time $t$.
\end{itemize}

\subsubsection*{Refined definition in case of full information}

In this subsection we focus on the special case of full information $\mathcal{G}=\mathcal{F}$.
State-wise prospective reserves $(V_j)_j$ are not necessarily unique. Theorem~\ref{Proposition:RegVj} only shows the uniqueness of  the processes $(I_j V_j)_j$. This is sufficient for the stochastic Thiele equations~\eqref{PRE:Thiele} and~\eqref{PRE_generalizedThiele}, but the representation~\eqref{PRE:SumAtRiskSimplified} entails the term $t \mapsto \mathbf{1}_{\{Z_{t-} = j\}} V^{\mathcal{F}}_k(t)$ for $k\neq j$. So for~\eqref{PRE:SumAtRiskSimplified} we need to specify the state-wise prospective reserves more carefully.

While the state-wise prospective reserves defined by \eqref{RegDefOfVj} do not necessarily satisfy~\eqref{PRE:SumAtRiskSimplified}, the following alternative construction meets that equation.
\begin{theorem}\label{Thm:VFjConstruction}
Let $\zeta_x:=\xi_x$, $x \in \mathcal{S}$, such that $\mathcal{F}=\mathcal{G}$. The processes $V_j=(V_j(t))_{t\geq0}$, $j\in S$, defined by
\begin{align}\label{VFjDefinition}
  V_j(t)&:=  \left\{ \begin{array}{ll}
  \sum_{n=0}^{\infty} \mathbf{1}_{\{J^j_t=n\}} \mathcal{V}_{(j,n)}(t) &\text{ if }  t>0, Z_{t-}=j \text{ or } t=0, \\
  V_{kj}(t) &\text{ if }  t>0, Z_{t-}=k \neq j, \end{array} \right.
\end{align}
are state-wise prospective reserves. Moreover,
\begin{align}\label{VjVjjVjk}\begin{split}
V_j(t)  &=I_j(t-)\, V_{jj}(t)+ \sum_{k:k\neq j} I_k(t-)\, V_{kj}(t), \quad j\in S,t>0,
\end{split}\end{align}
almost surely.
\end{theorem}
\begin{proof}
For $t>0$, $Z_{t-}= j$  the definition of $V_j(t)$ equals \eqref{RegDefOfVj}, and
 we inherit the results from Theorem \ref{Proposition:RegVj}. Based on definition~\eqref{VFjDefinition} and using Proposition~\ref{Prop:MeaningOfVjjVjkVkj}, we almost surely obtain that
\begin{align*}
 I_k(t-) I_j(t)\, V_j(t) = \sum_{k:k\neq j} I_k(t-) I_j(t)\, V_{kj}(t) = I_k(t-) I_j(t)\, V(t),
\end{align*}
so in combination with Theorem  \ref{Proposition:RegVj} we can conclude that $V_j$, $j\in S$, are indeed state-wise prospective reserves.

The simultaneous almost sure equation $I_j(t-) \,V_k(t)  =I_j(t-)\, V_{jk}(t)$, $t>0$ is a  direct consequence of definition~\eqref{VFjDefinition}.
Because of $\sigma(\xi_{x}) \vee \mathcal{A} = \mathcal{F}_{T_{x}}$ and $\mathcal{I}_{x}(t-) \,\mathcal{I}_{x}(t)= \mathbf{1}_{\{T_x < t\}} \,\mathcal{I}_{x}(t) $ for $x= (j,n)$,  we moreover can show that
\begin{align*}
  I_j(t-)\, V_{jj}(t) &=
  \sum_{n=0}^{\infty} \mathbf{1}_{\{J^j_{t-}=n\}} \mathcal{I}_{(j,n)}(t-) \mathcal{V}_{(j,n)(j,n)}(t)\\
  &=\sum_{n=0}^{\infty} \mathbf{1}_{\{J^j_{t-}=n\}} \mathcal{I}_{(j,n)}(t-) \mathcal{V}_{(j,n)}(t)\\
  &=I_j(t-)V_{j}(t)
\end{align*}
almost surely. By similar arguments as in the proof of Proposition~\ref{Proposition:RegVx}, property (3) of Definition~\ref{def:stat_counters} implies that the paths of $V_{j}$  are right-continuous on the intervals where $I_j$ equals one. Similarly, one can show this property also for $V_{jj}$. As a consequence, $ I_j(t-)\, V_{jj}(t)= I_j(t-)V_{j}(t)$ is not only point-wise almost surely true but also simultaneously almost surely true for all $t>0$.
\end{proof}
Note that \eqref{VjVjjVjk} and Proposition \ref{Prop:MeaningOfVjjVjkVkj} imply that
\begin{align*}
  V_j(t) = \E[ Y(t) \, | \, Z_{t-},F_{t-},Z_t=j]
\end{align*}
almost surely for each $t\geq0$, which can be intuitively expressed as
\begin{align*}
  V_j(t) = \E[ Y(t) \, | \, \mathcal{F}_{t-}, Z_{t}=j ]
\end{align*}
since $\mathcal{F}_{t-}=\sigma( Z_{t-},F_{t-}) \vee \mathcal{A}$.
The latter kind of expression is used in equations (3.5) and (4.6) of~\citet{ChristiansenDjehiche2020}, and we propose here to interpret it as the process $V_j$ according to Theorem~\ref{Thm:VFjConstruction}.

\subsubsection*{State-wise payment functions}

In the setting with information discarding, the $\mathcal{F}$-predictable
processes $b_j$ and $b_{jk}$ might be unobservable with respect to the reduced information $\mathcal{G}$. Thus we need projections of $b_j$ and $b_{jk}$ into the reduced information space $\mathcal{G}$. In the following, we first only provide an informal (point-wise almost sure) definition; the extension to the rigorous process point of view is then covered by the subsequent Proposition~\ref{JointlyMeasurableVersions2}. The resulting quantities are mathematically very similar to the previously investigated state-dependent quantities pertaining to the prospective reserve, so we keep this subsection relatively brief.

For $j,k\in S$, $j\neq k$, $g\in E$, we denote $a_j$ given by
\begin{align}\label{DefStateWisePayments11}
a_j(t) =      \E[  b_j(t) \, | \,  Z_{t-}=j,G^j_{t-}], \quad t\geq 0,
\end{align}
as the \emph{$\mathcal{G}$-averaged sojourn payment rate in state $j$}, and we denote $a_{jk}(\cdot,g)$ given by
\begin{align}\label{DefStateWisePayments12}
a_{jk}(t,g) = \E[  b_{jk}(t) \, | \, Z_{t-}=j,G^j_{t-},Z_{t}=k,G^k_{t}=g ], \quad t\geq0,
\end{align}
as the \emph{$\mathcal{G}$-averaged transition payment upon a transition from $j$ to $k$ with $G^k_t=g$}. 

These quantities are cast in the terms of the state process $Z$. Analogously, we can consider $\mathcal{G}$-averaged payments in terms of the expanded state process $\mathcal{Z}=(Z,J)$. For $x,y\in S$, $x\neq y$, $e\in E'$ with $x=(j,n)$ and $y=(k,l)$, these are given according to
\begin{align}\label{DefStateWisePayments21}
\mathpzc{a}_{\,(j,n)}(t) &=      \E[  b_j(t) \, | \,\mathcal{Z}_{t-}=(j,n),\zeta_{(j,n)}], \quad t\geq0, \\ \label{DefStateWisePayments22}
\mathpzc{a}_{\,(j,n)(k,l)}(t,e) &= \E[  b_{jk}(t) \,| \, \mathcal{Z}_{t-}=(j,n), \zeta_{(j,n)},\mathcal{Z}_{t}=(k,l),\zeta_{(k,l)}=e ], \quad t\geq0.
\end{align}

\begin{proposition}\label{JointlyMeasurableVersions2}
There exist modifications of the $\mathcal{G}$-average payments given via~\eqref{DefStateWisePayments11}--\eqref{DefStateWisePayments22} that are jointly measurable as mappings of $(\omega,t)$, $(\omega, t,g)$, $(\omega,t)$, and $(\omega, t,e)$, respectively.
\end{proposition}
\begin{proof}[Sketch of proof]
The desired expressions can be defined analogously to Definition~\ref{DefFurtherStateWiseQuant1} and Definition~\ref{DefFurtherStateWiseQuant2} by using the same arguments as in Proposition~\ref{Prop:MeaningOfVjjVjkVkj} and Proposition~\ref{Prop:MeaningOfVxxVxyVyx}. In particular, define $\mathpzc{a}_{\,(j,n)}$ according to
\begin{align*}
\mathpzc{a}_{\,(j,n)}(t) :=& \frac{\E[ \mathcal{I}_{(j,n)}(t-)b_j(t) \, | \, \zeta_{(j,n)}]}{\E[ \mathcal{I}_{(j,n)}(t-) \, | \, \zeta_{(j,n)}]}, \quad t\geq0.\qedhere
\end{align*}
\end{proof}
In the following, $a_j$, $a_{jk}$, $\mathpzc{a}_{\,(j,n)}$, $\mathpzc{a}_{\,(j,n)(k,l)}$ always refer to the modifications indicated in the proof of Proposition~\ref{JointlyMeasurableVersions2}.

\section{Infinitesimal compensators}\label{sec:inf_comps}

The classic stochastic Thiele equation can be derived by applying martingale techniques, in particular by using the concept of compensators, see e.g.~\citet{Norberg1992} and \citet{ChristiansenDjehiche2020}. Since martingale theory is generally restricted to monotone information only, in our non-monotone information setting we need to resort to the so-called infinitesimal approach suggested by~\citet{Christiansen2020}. In the infinitesimal approach, the classic concept of compensators is replaced by so-called infinitesimal forward/backward compensators.

The ordinary compensator w.r.t.~$\mathcal{F}$ of an adapted process $X=(X(t))_{t\geq0}$ with c\`{a}dl\`{a}g paths and finite expected variation on compacts is a process $C$ that satisfies
\begin{align*}
      C(\d t) = \E[ X(\d t) \, | \, \mathcal{F}_{t-}], \quad C(0)=0,
    \end{align*}
which is informal notation for the point-wise almost sure identity
\begin{align*}
      C(t) = \lim_{n\rightarrow \infty} \sum_{\mathfrak{T}_n} \E[ X(t_{k+1}) - X(t_{k}) \, | \, \mathcal{F}_{t_{k}} ],
\end{align*}
      where $(\mathfrak{T}_n)_{n\in \mathbb{N}}$ is any non-decreasing sequence (i.e.~$\mathfrak{T}_n \subseteq \mathfrak{T}_{n+1}$ for all $n$) of partitions $0=t_0 < \cdots < t_n=t$ of the interval $[0,t]$ such that $\lvert\mathfrak{T}_n\rvert := \max\{ t_{k}-t_{k-1}: k=1, \ldots, n\} \rightarrow 0$ for $n \rightarrow \infty$. Following Christiansen (2020), we denote $C$ and $\overline{C}$ defined point-wise by
    \begin{align*}
      C(\d t) &= \E[ X(\d t) \, | \, \mathcal{G}_{t-}],  \quad C(0)=0,\\
      \overline{C}(\d t) &= \E[ X(\d t) \, | \, \mathcal{G}_{t}],\quad \overline{C}(0)=0,
    \end{align*}
as the infinitesimal forward compensator and infinitesimal backward compensator of $X$ w.r.t.~$\mathcal{G}$, respectively. Again, the notation alludes to the following point-wise almost sure identities:
\begin{align*}
      C(t) &= \lim_{n\rightarrow \infty} \sum_{\mathfrak{T}_n} \E[ X(t_{k+1}) - X(t_{k}) \, | \, \mathcal{G}_{t_{k}} ], \\
      \overline{C}(t) &= \lim_{n\rightarrow \infty} \sum_{\mathfrak{T}_n} \E[ X(t_{k+1}) - X(t_{k}) \, | \, \mathcal{G}_{t_{k+1}} ].
\end{align*}
In the case of full information $\mathcal{G}=\mathcal{F}$, the infinitesimal forward compensator equals the ordinary compensator and the infinitesimal backward compensator equals the process itself, since $\E[ X(\d t) \, | \, \mathcal{F}_{t}]= X(\d t)$.

To derive the classic stochastic Thiele equation, one needs the compensators of the counting processes $N_{jk}$, $j,k \in S$, $j\neq k$. In our generalized setting, we need instead infinitesimal forward and backward compensators for relevant counting processes. For each $x,y\in \mathcal{S}$, $x\neq y$, $A' \in  \mathfrak{B}(E')$ let
\begin{align*}
 \mathcal{C}_{xy}(\d t \times A')&:= \mathcal{I}_x
(t-) \,\gamma_{xy}(\d t \times A'), \\
\overline{\mathcal{C}}_{xy}(\d t \times A')&:= \mathcal{I}_y
(t) \,\overline{\gamma}_{xy}(\d t \times A'),
\end{align*}
where the random measures $\gamma_{xy}$ and $\overline{\gamma}_{xy}$  are given as the unique completions of
\begin{align*}
\gamma_{xy} ((0,t] \times  A') &:= \!   \int_{(0,t] \times A'} \mathbf{1}_{\{\E[ \mathcal{I}_{x}(s-)  \, | \,  \zeta_x ] > 0\}} \frac{\E[ \mathcal{I}_{x}(s-)  \, | \,  \zeta_{x},T_y=s,\zeta_y=e]}{ \E[ \mathcal{I}_{x}(s-)  \, | \,   \zeta_{x} ]} \,  P( (T_{y},\zeta_y) \in \d s \times \d e \, | \, \zeta_x),\\
\overline{\gamma}_{xy} ((0,t] \times  A') &:= \!   \int_{(0,t] \times A'} \mathbf{1}_{\{\E[ \mathcal{I}_{y}(s) \, | \,  \zeta_y ] > 0\}}  \frac{\E[ \mathcal{I}_{y}(s)  \, | \,  \zeta_{y},S_x=s,\zeta_x=e]}{ \E[ \mathcal{I}_{y}(s) \, | \, \zeta_{y} ]} \, P( (S_{x},\zeta_x) \in \d s \times \d e \, | \, \zeta_y).
\end{align*}
\begin{proposition}\label{prop:compensators_of_gamma}
For each $x,y\in \mathcal{S}$, $x\neq y$, $A' \in  \mathfrak{B}(E')$ it point-wise almost surely holds that
\begin{align*}
\mathcal{C}_{xy}(\d t \times A')
&=
\E[\mathbf{1}_{\{\zeta_y \in A'\}}\,\mathcal{N}_{xy}(\d t) \, | \, \mathcal{G}_{t-}], \\
\overline{\mathcal{C}}_{xy}(\d t \times A')
&=
\E[\mathbf{1}_{\{\zeta_x \in A'\}}\,\mathcal{N}_{xy}(\d t) \, | \, \mathcal{G}_{t}].
\end{align*}
\end{proposition}
\begin{proof}
See Proposition 5.1 and Theorem 5.2 in~\citet{Christiansen2020}.
\end{proof}
For each $j,k \in S$, $j\neq k$, $A=\{n \} \times A' \in  \mathfrak{B}(E)$ let
\begin{align*}
 C_{jk}(\d t \times A)&:= I_j
(t-) \,\lambda_{jk}(\d t \times A), \\
\overline{C}_{jk}(\d t \times A)&:= I_k(t) \,\overline{\lambda}_{jk}(\d t \times A),
\end{align*}
where the random measures $\lambda_{jk}$ and $\overline{\lambda}_{jk}$ are given as
\begin{align}\label{eq:randommeasure_jk}
 \lambda_{jk}(\d t \times A) &:=
 \,\gamma_{(j,J^j_{t-})(k,n)}(\d t \times A'),\\ \nonumber
 \overline{\lambda}_{jk}(\d t \times A) &:=
 \,\overline{\gamma}_{(j,n)(k,J^k_t)}(\d t \times A').
\end{align}
The random measures depend on the specification of $\zeta_x$, $x\in\mathcal{S}$. We denote by $\lambda_{jk}^\mathcal{F}$, $j,k\in S$, $j \neq k$, the random measures corresponding to $\zeta_x = \xi_x$, $x\in\mathcal{S}$.
\begin{proposition}
For each $j,k \in S$, $j\neq k$, $A=\{n \} \times A' \in  \mathfrak{B}(E)$ it point-wise almost surely holds that
\begin{align*}
C_{jk}(\d t \times A)
&=
\E[\mathbf{1}_{\{G_t \in A\}}\,N_{jk}(\d t) \, | \, \mathcal{G}_{t-}], \\
\overline{C}_{jk}(\d t \times A)
&=
\E[\mathbf{1}_{\{G_{t-} \in A\}}\,N_{jk}(\d t) \, | \, \mathcal{G}_{t}].
\end{align*}
\end{proposition}
\begin{proof}
Apply the summations $ \sum_{l=0}^{\infty} \mathcal{I}_{(j,l)}(t-)= I_j(t-)$ and $ \sum_{l=0}^{\infty} \mathcal{I}_{(j,l)}(t)= I_j(t)$ on the equations in Proposition~\ref{prop:compensators_of_gamma}.
\end{proof}
The next proposition, which succeeds a technical condition, illustrates the roles of $\gamma_{xy},\lambda_{jk}$ and $ \overline{\gamma}_{xy},\overline{\lambda}_{jk} $ as forward and backward transition intensities. 
\begin{condition}\label{cond:majorant}
There exist non-negative c\`{a}gl\`{a}d processes $\mu^{\mathcal{F}}_{jk}(\cdot,E)$, $j,k \in S$, $j\neq k$, which are on each compact interval simultaneously dominated by a random variable with finite expectation and satisfy $\lambda^{\mathcal{F}}_{jk}(\d t \times E) = \mu^{\mathcal{F}}_{jk}(t,E) \, \d t$ for all $j,k \in S$, $j\neq k$.\democd
\end{condition}
\begin{proposition}\label{prop:smoothness_fixed_again}
Suppose Condition~\ref{cond:majorant} holds.
\begin{itemize}
\item Let $x,y \in\mathcal{S}$, $x\neq y$, $A'\in\mathfrak{B}(E')$. If there exists a non-negative c\`{a}gl\`{a}d process $\rho_{xy}(\cdot,A')$ satisfying $\gamma_{xy}(\d t \times A') = \rho_{xy}(t,A')\, \d t$, then almost surely simultaneously for all $t\geq0$,
\begin{align*}
\rho_{xy}(t+, A') &= \lim_{h\downarrow 0}\frac{1}{h} P\big(\mathcal{Z}_{t+h}=y ,\zeta_{y}\in \times A' \, \big| \, \mathcal{Z}_t=x,\zeta_{x}\big).
\end{align*}
\item Let $x,y \in\mathcal{S}$, $x\neq y$, $A'\in\mathfrak{B}(E')$. If there exists a non-negative c\`{a}dl\`{a}g process $\overline{\rho}_{xy}(\cdot,A')$ satisfying $\overline{\gamma}_{xy}(\d t \times A') = \overline{\rho}_{xy}(t,A')\, \d t$,
then almost surely simultaneously for all $t>0$,
\begin{align*}
  \overline{\rho}_{xy}(t-, A') &= \lim_{h\downarrow 0}\frac{1}{h} P\big(\mathcal{Z}_{t-h}=x ,\zeta_{x}\in A' \,  \big| \, \mathcal{Z}_t=y,\zeta_{y}\big).
\end{align*}
\item Let $j,k \in S$, $j \neq k$, $A\in\mathfrak{B}(E)$. If there exists a non-negative c\`{a}gl\`{a}d process $\mu_{jk}(\cdot,A)$ satisfying $\lambda_{jk}(\d t \times A) = \mu_{jk}(t,A)\, \d t$, then almost surely simultaneously for all $t\geq0$,
\begin{align*}
  \mu_{jk}(t+, A) &= \lim_{h\downarrow 0}\frac{1}{h} P\big(Z_{t+h}=k, G_{t+h}\in A \, \big| \, Z_{t}=j,G^j_t\big).
\end{align*}
\item Let $j,k \in S$, $j \neq k$, $A\in\mathfrak{B}(E)$. If there exists a non-negative c\`{a}dl\`{a}g process $\overline{\mu}_{jk}(\cdot,A)$ satisfying $\overline{\lambda}_{jk}(\d t \times A) = \overline{\mu}_{jk}(t,A)\, \d t$, then almost surely simultaneously for all $t>0$,
\begin{align*}
  \overline{\mu}_{jk}(t-, A) &= \lim_{h\downarrow 0}\frac{1}{h}P\big(Z_{t-h}=j, G_{t-h} \in A  \, \big| \, Z_{t}=k,G^k_t\big).
\end{align*}
\end{itemize}
\end{proposition}
\begin{proof}[Sketch of proof]
Condition~\ref{cond:majorant} ensures that the probability of two or more jumps occurring in a time interval of length $h$ vanishes faster than $h$, see also the proof of Lemma 3.3 in~\citet{Aalen1978}. It is then possible to show that
   \begin{align*}
   \lim_{h\downarrow 0}\frac{1}{h}  P\big(\mathcal{Z}_{t+h}=y ,\zeta_{y}\in A' \, \big| \, \mathcal{Z}_t=x,\zeta_{x}\big)
    &= \lim_{h\downarrow 0}\frac{1}{h} \E \bigg[ \int_{(t,t+h]} \mathbf{1}_{\{\zeta_{y}\in A'\}} \, \mathcal{N}_{xy}(\d s) \, \bigg| \, \mathcal{Z}_t=x,\zeta_{x}\bigg]\\
    &= \lim_{h\downarrow 0}\frac{1}{h}  \int_{(t,t+h]} \gamma_{xy}(\d s \times A')
  \end{align*}
  almost surely. The path properties of $\rho_{xy}(\cdot,A')$ ensure that the result holds not only point-wise almost surely but also simultaneously for all $t\geq0$. This proofs the first limit. The proofs of the other limits are similar.
\end{proof}
Consider the discounted payment process $\tilde{B}=(\tilde{B}(t))_{t\geq0}$ given by $\tilde{B}(0):=B(0)$ and $\tilde{B}(\d t):= v(t) \, B(\d t)$. Let 
\begin{align*}
C_{\tilde{B}}(\d t) :=  \sum_{x\in\mathcal{S}} \mathcal{I}_x(t-)\,v(t)  \,\mathpzc{a}_x(t) \, \beta(\d t) +\sum_{x,y\in \mathcal{S}\atop x\neq y} \int_{E'} v(t) \, \mathpzc{a}_{\,xy}(t,e)
   \,  \mathcal{I}_x(t-) \,\gamma_{xy} ( \d t \times \d e).
\end{align*}
Note that by the summation $\sum_{l=0}^{\infty} \mathcal{I}_{(j,l)}(t-)= I_j(t-)$,
\begin{align*}
C_{\tilde{B}}(\d t)
=
\sum_{j \in S} I_j(t-)\,v(t)\, a_j(t) \, \beta(\d t) +\sum_{j,k\in S\atop j\neq k} \int_{  E} v(t) \, a_{jk}(t,g)
   \, I_j(t-)\, \lambda_{jk} ( \d t \times \d g).
\end{align*}
\begin{proposition}
It point-wise almost surely holds that
\begin{align*}
C_{\tilde{B}}(\d t)
=
\E[\tilde{B}(\d t) \, | \, \mathcal{G}_{t-}].
\end{align*}
\end{proposition}
\begin{proof}See Theorem 5.2 and Example 7.2 in~\citet{Christiansen2020}.
\end{proof}

\section{Stochastic Thiele equation} \label{sec:stoch_thiele}

In this section, we present the main result of the paper by deriving a generalization of the stochastic Thiele equation to allow for non-monotone information. The derivation relies on the explicit infinitesimal martingale representation \citep[see Theorem 6.1 and Theorem 7.1 in][]{Christiansen2020}. In comparison, the classic stochastic Thiele equation with full information $\mathcal{G}=\mathcal{F}$ is closely related to the classic martingale representation theorem, see e.g.~\citet{Norberg1992} and \citet{ChristiansenDjehiche2020}.
\begin{theorem}\label{thm:sdeVx}
The processes $\mathcal{V}_x$, $x \in \mathcal{S}$, almost surely solve the stochastic differential equation
\begin{align*}
      0=\sum_{x\in \mathcal{S}} \mathcal{I}_x(t-) \bigg( \mathcal{V}_x(\d t) -  \mathcal{V}_x(t-) \frac{\kappa(\d t)}{\kappa(t-)} + \mathpzc{a}_x(t)\, \beta(\d t)
       &+  \sum_{y: y\neq x} \int_{ E'} \mathcal{R}_{xy}(t,e) \,\gamma_{xy}(\d t \times  \d e) \\
       &-  \sum_{y:y\neq x} \int_{E'} \overline{\mathcal{R}}_{yx}(t,e)\,   \overline{\gamma}_{yx}(\d t \times \d e) \bigg),
 \end{align*}
where for $x,y\in\mathcal{S}$, $x\neq y$, $e\in E'$, $t\geq0$,
    \begin{align*}
        \mathcal{R}_{xy}(t,e)&:= \mathpzc{a}_{\,xy}(t,e) + \mathcal{V}_{xy}(t,e) -\mathcal{V}_{xx}(t), \\
        \overline{\mathcal{R}}_{yx}(t,e)&:= \overline{\mathcal{V}}_{yx}(t,e) -\mathcal{V}_{xx}(t).
    \end{align*}
\end{theorem}
The processes  $\mathcal{R}_{xy}(\cdot,e)$ and $ \overline{\mathcal{R}}_{xy}(\cdot,e)$ may be interpreted as the forward sum at risk and backward sum at risk upon a transition of $\mathcal{Z}$ from $x$ to $y$ with updated and original information $e$, respectively.

Recall that all processes appearing in Theorem~\ref{thm:sdeVx} have been carefully defined in the previous sections. This is also the case for the processes appearing in Theorem~\ref{thm:sdeVj} and Corollary~\ref{cor:monotone_Thiele} below. By taking such great care, we have guaranteed their sufficient regularity, which ensures that the remainder of derivations in this paper, including the proofs of the aforementioned results, can proceed without any unease regarding uncountability of null-sets; in particular, whenever strictly necessary, intermediate results can be shown to hold not only almost surely in a point-wise but also in a simultaneous sense. With this settled, and with the intent of improving readability, we in general from this point onward cease to be explicit about null-sets.
\begin{theorem}[Generalized stochastic Thiele equation]\label{thm:sdeVj}
The processes $V_j$, $j \in S$, defined via~\eqref{RegDefOfVj} almost surely solve the stochastic differential equation
    \begin{align*}
      0=\sum_{j\in S} I_j(t-) \bigg( V_j(\d t) -  V_j(t-) \frac{\kappa(\d t)}{\kappa(t-)} + a_j(t)\, \beta(\d t)
       &+  \sum_{k: k\neq j} \int_{ E} R_{jk}(t,g) \,\lambda_{jk}(\d t \times  \d g) \\
       &-  \sum_{k:k\neq j} \int_{ E} \overline{R}_{kj}(t,g)\,   \overline{\lambda}_{kj}(\d t \times \d g) \bigg),
    \end{align*}
where for $j,k\in S$, $j\neq k$, $g\in E$, $t\geq0$,
    \begin{align*}
        R_{jk}(t,g)&:= a_{jk}(t,g) + V_{jk}(t,g)-V_{jj}(t), \\
        \overline{R}_{kj}(t,g)&:= \overline{V}_{kj}(t,g)-V_{jj}(t).
    \end{align*}
\end{theorem}
We may interpret $R_{jk}(t,g)$ and $ \overline{R}_{jk}(t,g)$ as the forward sum at risk and backward sum at risk upon a transition from $j$ to $k$ with updated and original information $g$, respectively.

For $j,k \in S$, $j\neq k$, let
\begin{align*}
\lambda_{jk}(\mathrm{d}t) := \lambda_{jk}(\mathrm{d}t \times E).
\end{align*}
\begin{corollary}[Classic stochastic Thiele equation]\label{cor:monotone_Thiele}
Let $\zeta_x=\xi_x$, $x \in \mathcal{S}$, such that $\mathcal{G}= \mathcal{F}$. Then the processes $V_j$, $j \in S$, defined via~\eqref{VFjDefinition} almost surely solve the stochastic differential equation
\begin{align*}
      0=\sum_{j\in S} I_j(t-) \bigg( V_j(\d t) -  V_j(t-) \frac{\kappa(\d t)}{\kappa(t-)} + b_{j}(t)\, \beta(\d t)
       &+ \sum_{k: k\neq j} R_{jk}(t)\, \lambda_{jk}(\d t ) \bigg),
    \end{align*}
where for $j,k\in S$, $j \neq k$, $t\geq0$,
\begin{align*}
   R_{jk}(t)&:=b_{jk}(t) + V_{k}(t) -V_{j}(t).
\end{align*}
\end{corollary}

Before we present the proofs of Theorem~\ref{thm:sdeVx}, Theorem~\ref{thm:sdeVj}, and Corollary~\ref{cor:monotone_Thiele}, we first provide an interpretation of the results. In the presence of full information $\mathcal{G}=\mathcal{F}$, Corollary~\ref{cor:monotone_Thiele} yields a stochastic differential equation that is directly comparable to the stochastic Thiele equation from~\citet{Norberg1992,Norberg1996}. In~\citet{Norberg1992,Norberg1996}, the compensators of the multivariate counting process are assumed to admit densities w.r.t.\ the Lebesgue-measure, and the result is derived by suitably applying the martingale representation theorem and identifying the integrands. The method of the present paper, while extended to also cover the case of non-monotone case, is based on a suitable application of the explicit infinitesimal martingale representation theorem. Actually, Corollary~\ref{cor:monotone_Thiele} may also be derived directly from the classic martingale representation theorem following~\citet{ChristiansenDjehiche2020}; in this case, the restriction to payments on the form~\eqref{eq:nice_payments} is not necessary.

The stochastic differential equations of Theorem~\ref{thm:sdeVx} and Theorem~\ref{thm:sdeVj} are in a twofold manner fundamentally different from the classic stochastic Thiele equation of Corollary~\ref{cor:monotone_Thiele}. Firstly, the sum at risks appearing in the term involving infinitesimal forward compensators, which correspond to ordinary compensators in the presence of full information, take a different form. Rather than being the difference of two state-wise prospective reserves added the relevant transition payment, it involves the difference of the forward state-wise prospective reserve and the prospective reserve for staying in the state added the relevant transition payment. In the presence of monotone information, we can show that the forward state-wise prospective reserve and the prospective reserve for staying in the state can be replaced by relevant ordinary state-wise prospective reserves, but this is not necessarily the case in the presence of non-monotone information. Here the possibility of information discarding entails a possible improvement in the accuracy of the reserving by utilizing the information available at time $t-$ and time $t$, rather than utilizing only the information available at time $t$.

Secondly, the stochastic differential equations of Theorem~\ref{thm:sdeVx} and Theorem~\ref{thm:sdeVj} contain an additional term that relates to the infinitesimal backward compensators. In the presence of monotone information, we can show that this term is zero. It is the backward looking equivalent of the term involving the infinitesimal forward compensators. Based on the information currently available, the term adjusts the dynamics to take into account the possibility that information discarding has just occurred.

In Section~\ref{sec:examples} we derive and interpret generalized stochastic Thiele equations in the presence of specific examples of non-monotone information related to stochastic retirement. We refer to this section for further interpretation and discussion of the general results.

\begin{proof}[Proof of Theorem~\ref{thm:sdeVx}]
According to Theorem 7.1 in~\citet{Christiansen2020}, the process $\tilde{V} := v \cdot V$ satisfies the equation
\begin{align*}
\tilde{V}(\d t)= - C_{\tilde{B}}(\d t)
&+ \sum_{x,y \in \mathcal{S}\atop x\neq y} \int_{E'} v(t) \big(  \mathcal{V}_{xy}(t,e) - \mathcal{V}_{xx}(t)\big) (\nu_{xy}-C_{xy})(\d t \times \d e) \\
&  -\sum_{x,y \in \mathcal{S} \atop x\neq y} \int_{E'} v(t) \big(  \overline{\mathcal{V}}_{xy}(t,e) - \mathcal{V}_{xx}(t) \big) (\overline{\nu}_{xy}-\overline{C}_{xy})(\d t \times \d e)
\end{align*}
for random measures $\nu_{xy}$ and $\overline{\nu}_{xy}$  that are defined as the unique completions of
\begin{align*}
\nu_{xy} ( (0,t] \times  A) &=\int_{(0,t]} \mathbf{1}_{\{\zeta_y \in A\}} \, \mathcal{N}_{xy}(\d s),\\
  \overline{\nu}_{xy} ( (0,t] \times  A) &=\int_{(0,t]} \mathbf{1}_{\{\zeta_x \in A\}} \, \mathcal{N}_{xy}(\d s).
\end{align*}
On the other hand, by applying integration by parts on $\tilde{V}(t) = \sum_{x\in\mathcal{S}} \mathcal{I}_x(t) \tilde{\mathcal{V}}_x(t)$, where $\tilde{\mathcal{V}}_x := v \cdot \mathcal{V}_x$, and using $\tilde{V} = v \cdot V$, we can show that
\begin{align*}
\tilde{V}(\d t) = \sum_{x\in\mathcal{S}} \mathcal{I}_x(t-) \tilde{\mathcal{V}}_x(\d t) + \sum_{x,y \in \mathcal{S} \atop x \neq y} v(t) \big( \mathcal{V}_y(t) - \mathcal{V}_x(t)\big) \nu_{xy}(\d t \times E').
\end{align*}
By equating the latter two equations and rearranging the terms, while using the fact that $C_{xy}(\d t \times \d e) = \mathcal{I}_x(t-)\,C_{xy}(\d t \times \d e)$ and $\overline{C}_{yx}(\d t \times \d e) =\mathcal{I}_x(t)\, \overline{C}_{yx}(\d t \times \d e)$ and the equation $\mathcal{I}_x(t) =  \mathcal{I}_x(t-) \,\mathcal{I}_x(t) + \mathbf{1}_{\{\mathcal{Z}_{t-}\neq x\}}\, \mathcal{I}_x(t)$,  we obtain
\begin{align*}
0
=&
\sum_{x\in\mathcal{S}} \mathcal{I}_x(t-) \Bigg(\tilde{\mathcal{V}}_x(\d t)  + v(t) \, \mathpzc{a}_x(t) \, \beta (\d t)
+ \sum_{y: y\neq x}\int_{E'} v(t) \, \mathcal{R}_{xy}(t,e) \, C_{xy}(\d t \times \d e) \\
&\qquad \qquad \qquad   -\sum_{y:y\neq x}\int_{E'} v(t) \, \overline{\mathcal{R}}_{yx}(t,e) \, \overline{C}_{yx}(\d t \times \d e) \Bigg) \\
& - \sum_{x,y \in \mathcal{S} \atop x \neq y }\int_{E'} v(t) \, \mathcal{V}_{xy}(t,e)  \,\nu_{xy}(\d t \times \d e) +  \sum_{x,y \in \mathcal{S} \atop x \neq y } \int_{E'} v(t) \,\overline{\mathcal{V}}_{xy}(t,e) \, \overline{\nu}_{xy}(\d t \times \d e)\\
& + \sum_{x,y \in \mathcal{S} \atop x \neq y } v(t)\big(   \mathcal{V} _{xx}(t)- \mathcal{V} _{yy}(t) \big) \nu_{xy}(\d t \times E')\\
&  -\sum_{x,y \in \mathcal{S} \atop y\neq x}\int_{E'} v(t) \big(  \overline{\mathcal{V}}_{yx}(t,e) - \mathcal{V} _{xx}(t) \big) \mathbf{1}_{\{\mathcal{Z}_{t-}\neq x \}} \, \mathcal{I}_x(t) \, \overline{C}_{yx}(\d t \times \d e)\\
& + \sum_{x,y \in \mathcal{S} \atop x\neq y}  v(t)\big( \mathcal{V}_{y}(t) - \mathcal{V}_{x}(t) \big)\nu_{xy}(\d t \times E').
\end{align*}
The third line equals $- v(t) \, V(t) \, \dot{N}(\d t) + v(t) \, V(t) \, \dot{N}(\d t) =0$ for $\dot{N} := \sum_{j,k\in S, j\neq k} N_{jk}$. By applying Proposition~\ref{prop:compensators_of_gamma} and the identity
\begin{align*}
\overline{\gamma}_{yx}(\d t \times \d e) = \mathbf{1}_{\{\mathcal{Z}_{t}\neq x\}} \overline{\gamma}_{yx}(\d t \times \d e) + \overline{C}_{yx}(\d t \times \d e),
\end{align*}
we further obtain
\begin{align*}
0
=&
\sum_{x\in\mathcal{S}} \mathcal{I}_x(t-) \Bigg(\tilde{\mathcal{V}}_x(\d t)  + v(t)\,\mathpzc{a}_x(t) \, \beta (\d t)
+ \sum_{y: y\neq x}\int_{E'} v(t)\, \mathcal{R}_{xy}(t,e) \, \gamma_{xy}(\d t \times \d e) \\
&\qquad \qquad \qquad   -\sum_{y:y\neq x}\int_{E'} v(t)\,\overline{\mathcal{R}}_{yx}(t,e) \, \overline{\gamma}_{yx}(\d t \times \d e) \Bigg) \\
& + \sum_{x,y \in \mathcal{S} \atop x \neq y } v(t)\big(   \mathcal{V} _{xx}(t)- \mathcal{V} _{yy}(t) \big) \nu_{xy}(\d t \times E')\\
&  -\sum_{x,y \in \mathcal{S} \atop y\neq x}\int_{E'} v(t) \big(  \overline{\mathcal{V}}_{yx}(t,e) - \mathcal{V} _{xx}(t) \big) \mathbf{1}_{\{\mathcal{Z}_{t-}\neq x \}}\,\mathcal{I}_x(t) \, \overline{\gamma}_{yx}(\d t \times \d e)\\
& + \sum_{x,y \in \mathcal{S} \atop x\neq y}  v(t)\big( \mathcal{V}_{y}(t) - \mathcal{V} _{x}(t) \big)\nu_{xy}(\d t \times E') \\
& +
\sum_{x,y \in \mathcal{S} \atop x \neq y} \int_{E'} v(t)\,\overline{\mathcal{R}}_{yx}(t,e)\, \mathcal{I}_x(t-)  \mathbf{1}_{\{\mathcal{Z}_{t}\neq x\}} \overline{\gamma}_{yx}(\d t \times \d e).
\end{align*}
 The third, forth, fifth and sixth line together equal
\begin{align*}
&  \sum_{x,y \in \mathcal{S}\atop x\neq y } v(t) (  \mathcal{V} _{y}(t) - \mathcal{V} _{yy}(t)) \, \nu_{xy}(\d t \times E')  \\
&  -\sum_{x,y \in \mathcal{S}\atop y\neq x}\bigg(\int_{E'} v(t)\big(  \overline{\mathcal{V}}_{yx}(t,e) - \mathcal{V} _{xx}(t) \big)\bigg(\sum_{z:z \neq x } \big(\nu_{zx}(\{ t\} \times E') - \nu_{xz}(\{t\} \times E')\big)\bigg) \overline{\gamma}_{yx}(\d t \times \d e)  \bigg)\\
&  +\sum_{x,y \in \mathcal{S}\atop x \neq y } v(t) (  \mathcal{V} _{xx}(t) - \mathcal{V} _{x}(t)) \, \nu_{xy}(\d t \times E')  \\
=&  \sum_{x,y \in \mathcal{S}\atop x\neq y } v(t) (  \mathcal{V} _{x}(t) - \mathcal{V} _{xx}(t)) \, \nu_{yx}(\d t \times E')  \\
&   -\sum_{x,y \in \mathcal{S}\atop y\neq x}\bigg( \sum_{z:z \neq x } \int_{E'} v(t)\big( \overline{\mathcal{V}}_{zx}(t,e) - \mathcal{V} _{xx}(t) \big)\, \overline{\gamma}_{zx}(\{t\} \times \d e)   \bigg)\big(\nu_{yx}(\d t \times E') - \nu_{xy}(\d t \times E')\big)\\
&   +\sum_{x,y \in \mathcal{S}\atop x \neq y } v(t) (  \mathcal{V} _{xx}(t) - \mathcal{V} _{x}(t)) \, \nu_{xy}(\d t \times E'),
\end{align*}
because $\sum_{z:z \neq x } \nu_{zx}(\{ t\} \times E')$ and $\sum_{z:z \neq x } \nu_{xz}(\{ t\} \times E')$ are non-zero only at finitely many time points. The latter three lines also add up to zero since
\begin{align*}
     \mathcal{V} _{x}(t)
     &= \mathcal{V} _{xx}(t) \bigg( 1 - \sum_{z: z\neq x}  \overline{\gamma}_{zx}(\{t \} \times E')\bigg)+ \sum_{z: z\neq x} \int_{E'}   \overline{\mathcal{V}}_{zx}(t,e) \,   \overline{\gamma}_{zx}(\{t \} \times \d e).
\end{align*}
This identity is a consequence of the following observations. If $\E[ \mathcal{I}_x(t) \, | \, \zeta_{x} ] = 0$, then by definition, $\mathcal{V} _{x}(t)=0$, $\mathcal{V} _{xx}(t) = 0$, and $ \overline{\gamma}_{zx}(\{t\}\times \d e) = 0$ and the identity simply reads $0=0$. On the other hand, if $\E[ \mathcal{I}_x(t) \, | \, \zeta_{x} ] >0$, then
\begin{align*}
     \mathcal{V} _{x}(t)  &= \frac{\E[ \mathcal{I}_x(t-)\mathcal{I}_x(t) Y(t) \, | \, \zeta_x]}{\E[ \mathcal{I}_x(t) \, | \, \zeta_x]}+ \frac{\E[  \sum_{z:z\neq x}\mathcal{I}_z(t-)\mathcal{I}_x(t) Y(t) \, | \, \zeta_x]}{\E[ \mathcal{I}_x(t) \, | \, \zeta_x]}\\
     & = \frac{\E[ \mathcal{I}_x(t-)\mathcal{I}_x(t) Y(t) \, | \, \zeta_x]}{\E[ \mathcal{I}_x(t-)\mathcal{I}_x(t) \, | \, \zeta_x]}\frac{\E[ \mathcal{I}_x(t-)\mathcal{I}_x(t) \, | \, \zeta_x]}{\E[ \mathcal{I}_x(t) \, | \, \zeta_x]}\\
     & \quad +  \sum_{z: z\neq x}  \int_{E'}\frac{\E[ \mathcal{I}_x(t) Y(t) \, | \, (S_z,\zeta_z,\zeta_x)=(t,e,\zeta_x)]}{\E[ \mathcal{I}_x(t) \, | \, (S_z,\zeta_z,\zeta_x)=(t,e,\zeta_x)]} \, \overline{\gamma}_{zx}(\{t\} \times \d e) \\
     &= \mathcal{V} _{xx}(t) \bigg( 1 - \sum_{z: z\neq x}  \overline{\gamma}_{zx}(\{t \} \times E')\bigg)+ \sum_{z: z\neq x} \int_{E'}  \overline{\mathcal{V}}_{zx}(t,e) \,   \overline{\gamma}_{zx}(\{t \} \times \d e).
\end{align*}
All in all, we have
\begin{align*}
0
=&
\sum_{x\in\mathcal{S}} \mathcal{I}_x(t-) \Bigg(\tilde{\mathcal{V}}_x(\d t)  + v(t)\,\mathpzc{a}_x(t)\, \beta (\d t)
+ \sum_{y: y\neq x}\int_{E'} v(t)\, \mathcal{R}_{xy}(t,e) \, \gamma_{xy}(\d t \times \d e) \\
&\qquad \qquad \qquad   -\sum_{y:y\neq x}\int_{E'} v(t)\,\overline{\mathcal{R}}_{yx}(t,e)\, \overline{\gamma}_{yx}(\d t \times \d e) \Bigg).
\end{align*}
Now apply integration by parts on $\tilde{\mathcal{V}} _x=v \cdot \mathcal{V} _x$ and rearrange the terms in order to end up with the statement of the theorem.
\end{proof}

\begin{proof}[Proof of Theorem~\ref{thm:sdeVj}]
The result follows from Theorem~\ref{thm:sdeVx} by  the summation $ \sum_{l=0}^{\infty} \mathcal{I}_{(j,l)}(t-)= I_j(t-)$ and the construction principle~\eqref{RegDefOfVj}, which applies not only to the state-wise prospective reserves but to all involved quantities. In the summation of the second line of the stochastic differential equation we may replace $J^j_{t-}$ by $J^j_{t}$ due to the left-continuity of $J^j$.
\end{proof}

\begin{proof}[Proof of Corollary~\ref{cor:monotone_Thiele}]
Since the definitions of $V_j(t)$ in \eqref{RegDefOfVj} and
\eqref{VFjDefinition} are identical on $\{Z_{t-}=j\}$,  the stochastic
differential equation from Theorem \ref{thm:sdeVj} applies here.

As $\{T_x=S_y\} \in \sigma(\zeta_x)$ and  $\mathbf{1}_{\{T_x=S_y \}}
\mathcal{I}_{x}(S_y)=\mathcal{I}_{x}(S_y)$, we have
\begin{align*}
\E[ \mathcal{I}_{x}(t) \,  | \,  \xi_{x},S_y=t,\xi_x=e]
&=
\mathbf{1}_{\{S_y = T_x\}} \, \E[ \mathcal{I}_{x}(t) \,  | \,
\xi_{x},S_y=t,\xi_x=e].
\end{align*}
Moreover, since  $\{t> T_x=S_y\} \in \sigma(\zeta_x)$ and
$\mathcal{I}_{x}(t-) \mathbf{1}_{\{t > T_x=S_y  \}} =\mathcal{I}_{x}(t-) 
\mathbf{1}_{\{S_y = T_x\}} $, we obtain
\begin{align*}
\mathcal{I}_x(t-) \, \mathbf{1}_{\{T_x= S_y \}} \,  P( (S_{y},\xi_y) \in
\d t \times \d e \, | \, \xi_x) &= 0.
\end{align*}
All in all, this means that  $\mathcal{I}_x(t-) \,
\overline{\gamma}_{yx}(\d t \times \d e) = 0$, using
the fact that $\zeta_x = \xi_x$.
By the summation $ \sum_{l=0}^{\infty} \mathcal{I}_{(j,l)}(t-)= I_j(t-)$
according to the proof of Theorem~\ref{thm:sdeVj}, we moreover get
$I_j(t-) \overline{\lambda}_{kj}(\d t \times \d g) = 0$. Thus, the
second line of the stochastic differential equation in Theorem
\ref{thm:sdeVj} is zero.

When we update the information from $(Z_{t-},F_{t-})$ to $(Z_{t},F_{t})$
at time $t$, the part $F_t$ is redundant since
\begin{align}\label{RedundantInformation}
   \sigma( Z_{t}, F_{t}) = \mathcal{F}_t= \mathcal{F}_{t-} \vee \sigma
(Z_t) = \sigma( Z_{t-}, F_{t-}) \vee \sigma(Z_t).
\end{align}
As a result, $\lambda_{jk}(\d t \times  \d g)$ equals $\lambda_{jk}(\d t
)$ if and only if $g=F_t$ and zero else, and we have $I_j(t-)
R_{jk}(t,F_t) = I_j(t-) R_{jk}(t)$ according to
Theorem~\ref{Thm:VFjConstruction}. Consequently,
\begin{align*}
\int_{ E} R_{jk}(t,g) \,\lambda_{jk}(\d t \times  \d g) &=\int_{ E}
R_{jk}(t) \,\lambda_{jk}(\d t ).\qedhere
\end{align*}
\end{proof}
In case the payments $B$ themselves depend on the prospective reserve $V$, the (stochastic) Thiele equations rather than~\eqref{eq:classic_res_def} might serve as definition for the prospective reserve $V$, see e.g.\ \citet{DjehicheLofdahl2016} and \citet{ChristiansenDjehiche2020}. In the presence of full information, this point of view is encapsulated by the following result.
\begin{proposition}\label{Prop:true_Thiele_equations}
Let $\zeta_x=\xi_x$, $x \in \mathcal{S}$, such that $\mathcal{G}= \mathcal{F}$.
Let there be a maximal contract time $n<\infty$ such that $b_j$ and $b_{jk}$ are constantly zero on the interval $(n, \infty)$.
Suppose that $W_j=(W_j(t))_{t\geq0}$, $j \in S$, are $\mathcal{F}$-predictable bounded processes for which $t \mapsto \int_{(0,t]} I_{j}(s-) \, W_j(\mathrm{d}s)$ defines finite variation processes. If $W_j$, $j \in S$, almost surely satisfy the stochastic differential equations
\begin{align}\label{eq:ThieleBSDE}
 0 = I_j(t-)\bigg(\!
 W_j(\d t) - W_j(t-) \frac{\kappa(\d t)}{\kappa(t-)} +  b_j(t) \, \beta(\d t) + \hspace{-1mm}  \sum_{k : k \neq j} (b_{jk}(t) \! + \! W_k(t) \! - \! W_j(t)) \, \lambda_{jk}(\d t) \!\bigg)
\end{align}
with terminal conditions $W_j(n)=0$, $j \in S$, then they are state-wise prospective reserves in the sense of  Definition \ref{def:stat_counters}.
\end{proposition}
\begin{proof}
The assumptions immediately yield conditions (2) and (3) of Definition~\ref{def:stat_counters}. By applying integration by parts and the stochastic differential equations for the processes $W_j$, $j\in S$, we obtain
\begin{align*}
  &\d \big( v(t) W_{Z_t}(t) \big) \\
  &=   \sum_{j \in S} I_j(t-) \bigg(v(t) W_j(\d t)- W_j(t-) v(t)\frac{\kappa(\d t)}{\kappa(t-)} \bigg) +\sum_{j,k \in S \atop j\neq k} v(t)( W_k(t) - W_j(t)) \, N_{jk}(\d t)\\
   &=   - v(t) \, B(\d t) + \sum_{j,k \in S \atop j \neq k}   v(t)( b_{jk}(t)+ W_k(t) - W_j(t)) (N_{jk}-\lambda_{jk})(\d t).
\end{align*}
Since each process $t \mapsto b_{jk}(t)+ W_k(t) - W_j(t)$ is $\mathcal{F}$-predictable and bounded, the last term is an $\mathcal{F}$-martingale. Consequently, we obtain
\begin{align*}
v(t)W_{Z_t}(t)
   &= \E\bigg[v(t)\sum_{j \in S} I_j(t-) W_j(t) \, \bigg| \, \mathcal{F}_t\bigg]\\
  &=\E\bigg[v(t)\int_{(t,n]} \frac{\kappa(t)}{\kappa(s)} \, B(\d s) \, \bigg| \, \mathcal{F}_t\bigg]
  = v(t) V(t).
\end{align*}
Noting $v>0$ yields condition (1) of Definition~\ref{def:stat_counters}, namely $I_j(t)W_j(t)=I_j(t)V(t)$, completing the proof.
\end{proof}

\section{Stochastic retirement: Stochastic Thiele equations} \label{sec:examples}

In this section, we adopt the setup concerning stochastic retirement of Section~\ref{sec:setup_stoc_re} and consider the non-monotone information given by $\mathcal{G}^{(1)}$ and $\mathcal{G}^{(2)}$. For comparison, we also look at the full information $\mathcal{G}^{(0)}:=\mathcal{F}$. Recall that $\mathcal{G}^{(1)}$ corresponds to the case where upon retirement or death the insurer does not have access to or discards all previous records regarding the health status of the insured, while $\mathcal{G}^{(2)}$ even keeps no record of the time of retirement. Our aim is to compare the prospective reserves $V^{(i)}$, $i\in\{0,1,2\}$, by studying the (generalized) stochastic Thiele equations in the presence of (non-monotone) information $\mathcal{G}^{(0)}$, $\mathcal{G}^{(1)}$, and $\mathcal{G}^{(2)}$. In the following, the superscripts $(0)$, $(1)$, and $(2)$ always refer to the setting with information $\mathcal{G}^{(0)}=\mathcal{F}$, $\mathcal{G}^{(1)}$, and $\mathcal{G}^{(2)}$, respectively.

Denote by $V^{(0)}_j$, $j \in S$, the processes defined in~\eqref{VFjDefinition}. Since information discarding does not occur before retirement or death, we may for $i\in\{1,2\}$ set $V^{(i)}_j := V^{(0)}_j$, $j\in\{1,\ldots,\varsigma\}$, and indeed obtain that $(V^{(i)}_j)_{t\geq0,j\in S}$ are state-wise prospective reserves in the sense of Definition~\ref{def:stat_counters} (under information $\mathcal{G}^{(i)}$) with $V^{(i)}_j$, $j \in\{\text{p},\text{d}\}$, defined according to~\eqref{RegDefOfVj}. Similarly, $\lambda_{jk}^{(i)} = \lambda_{jk}^{(0)}$ for $i\in\{1,2\}$, $j\in\{1,\ldots,\varsigma\}$, $k\in S$, $k \neq j$.

For $i\in\{0,1,2\}$, $j\in\{1,\ldots,\varsigma\}$, Corollary~\ref{cor:monotone_Thiele} yields
\begin{align*}
0
=
I_j(t-)\Bigg(V_j^{(i)}(\d t) - V_j^{(i)}(t-) \frac{\kappa(\d t)
}
{\kappa(t-)
} + b_j(t) \, \beta(\d t) &+ \hspace{-2mm}  \sum_{k : \varsigma \geq k \neq j} \big(b_{jk}(t) + V^{(i)}_k(t) - V_j^{(i)}(t)\big) \lambda_{jk}^{(i)}(\d t)    \\
&+  \hspace{-2mm}  \sum_{k \in \{\text{p},\text{d}\}} \big(b_{jk}(t) + V^{(0)}_k(t) - V_j^{(i)}(t)\big) \lambda_{jk}^{(i)}(\d t)    \Bigg)
\end{align*}
For $i\in\{1,2\}$ the sum at risks for $k\in\{\text{p},\text{d}\}$ take an unusual form as they involve $V_{\text{p}}^{(0)}$ and $V_{\text{d}}^{(0)}$ rather than $V_{\text{p}}^{(i)}$ and $V_{\text{d}}^{(i)}$; this constitutes a fundamental difference that arises due to non-monotonicity of information. Intuitively, it just reflects utilization of all available information before retirement or death. The situation, on the other hand, becomes more involved after retirement or death, as the three following theorems serve to illuminate.

\begin{theorem}\label{cor:sdeY_G0}
The state-wise prospective reserves $V^{(0)}_{\emph{p}}$ and $V^{(0)}_{\emph{d}}$ almost surely satisfy
\begin{align*}
&0 = I_{\emph{p}}(t-) \bigg(V_{\emph{p}}^{(0)}(\d t) - V_{\emph{p}}^{(0)}(t-) \frac{\kappa(\d t)}
{\kappa(t-)} + b_{\emph{p}}(t) \, \beta(\d t) + \big(b_{\emph{p}\emph{d}}(t) + V_{\emph{d}}^{(0)}(t) - V_{\emph{p}}^{(0)}(t)\big) \lambda_{\emph{p}\emph{d}}^{(0)}(\d t) \bigg), \\
&0 = I_{\emph{d}}(t-) \bigg(V_{\emph{d}}^{(0)}(\d t) - V_{\emph{d}}^{(0)}(t-) \frac{\kappa(\d t)}
{\kappa(t-)} + b_{\emph{d}}(t) \, \beta(\d t)\bigg).
\end{align*}
\end{theorem}
\begin{proof}
Apply Corollary~\ref{cor:monotone_Thiele}.
\end{proof}
\begin{theorem}\label{cor:sdeY_G1}
The state-wise prospective reserves $V^{(1)}_{\emph{p}}$ and $V^{(1)}_{\emph{d}}$ almost surely satisfy
\begin{align*}
&0 = I_{\emph{p}}(t-) \bigg(V_{\emph{p}}^{(1)}(\d t) - V_{\emph{p}}^{(1)}(t-) \frac{\kappa(\d t)}
{\kappa(t-)} + a_{\emph{p}}^{(1)}(t) \, \beta(\d t) + \big(a_{\emph{p}\emph{d}}^{(1)}(t) + V_{\emph{p}\emph{d}}^{(1)}(t) - V_{\emph{p}\emph{p}}^{(1)}(t)\big) \lambda_{\emph{p}\emph{d}}^{(1)}(\d t) \bigg), \\
&0 = I_{\emph{d}}(t-) \bigg(V_{\emph{d}}^{(1)}(\d t) - V_{\emph{d}}^{(1)}(t-) \frac{\kappa(\d t)}
{\kappa(t-)} + a_{\emph{d}}^{(1)}(t) \, \beta(\d t)\bigg).
\end{align*}
\end{theorem}
\begin{theorem}\label{cor:sdeY_G2}
The state-wise prospective reserves $V^{(2)}_{\emph{p}}$ and $V^{(2)}_{\emph{d}}$ almost surely satisfy
\begin{align*}
&0 = I_{\emph{p}}(t-) \bigg(V_{\emph{p}}^{(2)}(\d t) - V_{\emph{p}}^{(2)}(t-) \frac{\kappa(\d t)}
{\kappa(t-)} + a_{\emph{p}}^{(2)}(t) \, \beta(\d t) + \big(a_{\emph{p}\emph{d}}^{(2)}(t) + V_{\emph{p}\emph{d}}^{(2)}(t) - V_{\emph{p}\emph{p}}^{(2)}(t)\big) \lambda_{\emph{p}\emph{d}}^{(2)}(\d t) \\
 &\hspace{89mm}- \sum_{k \leq \varsigma}\int_E \overline{R}_{k\emph{p}}^{(2)}(t,g) \, \overline{\lambda}_{k\emph{p}}^{(2)}(\d t \times \d g)\bigg), \\
&0 = I_{\emph{d}}(t-) \bigg(V_{\emph{d}}^{(2)}(\d t) - V_{\emph{d}}^{(2)}(t-) \frac{\kappa(\d t)}
{\kappa(t-)} + a_{\emph{d}}^{(2)}(t) \, \beta(\d t)\bigg).
\end{align*}
\end{theorem}
\begin{proof}[Sketch of proof of Theorem~\ref{cor:sdeY_G1} and Theorem~\ref{cor:sdeY_G2}]
Apply Theorem~\ref{thm:sdeVj} and calculate some of the terms more explicitly.
\end{proof}

Recall that $\mathcal{G}^{(2)}$ does not have the time since retirement as admissible information. Such substantial information discarding leads to an additional term in the stochastic differential equation for $V_{\text{p}}^{(2)}$ related to the infinitesimal backward compensators $\overline{C}_{k\text{p}}^{(2)}(\d t \times \d g) = I_{\text{p}}(t) \, \overline{\lambda}_{k\text{p}}^{(2)}(\d t \times \d g)$, $k\in\{1,\ldots,\varsigma\}$. Recall that the time of retirement and death are given by the first hitting times $\eta$ and $\delta$, respectively. Straightforward calculations yield the following identity:
\begin{align}\begin{split}\label{remark:backward_sum_at_risk}
&\sum_{k \leq \varsigma}\int_E \overline{R}_{k\text{p}}^{(2)}(t,g) \, \overline{\lambda}_{k\text{p}}^{(2)}(\d t \times \d g) \\
&=
\Big(\E[ Y(t) \, | \, \eta=t] -  \mathbf{1}_{\{P(\eta < t < \delta )>0\}}\E[ Y(t) \, | \, \eta < t < \delta]\Big)
\frac{\mathbf{1}_{\{P(\eta \leq t < \delta )>0\}}}{P(\eta \leq t < \delta )} P(\eta\in\d t).
\end{split}
\end{align}
So the term adjusts the dynamics to take into account the possibility that retirement might just have occurred rather than having occurred some time ago (conditionally on the insured presently being retired). In the former case, at time $t$ one would reserve $\E[ Y(t) \, | \, \eta=t]$, while in the latter case one would reserve $E[ Y(t) \, | \, \eta < t < \delta]$. This constitutes a description of the first part of the product. The second part is exactly the infinitesimal probability of retirement having just occurred, conditionally on the insured presently being retired.

\section{Stochastic retirement: Feynman-Kac formulas} \label{sec:examples2}

We now specialize and simplify the setting to provide a more straightforward and less technical discussion of the general results and their relation to actuarial practice. We start out by specifying the model under full information $\mathcal{G}=\mathcal{G}^{(0)} = \mathcal{F}$. First, we suppose that $b_j$ and $b_{jk}$ are deterministic for all $j,k \in S$, $j\neq k$, and we assume the existence of a maximal contract time $n<\infty$ such that $b_j$ and $b_{jk}$ are constantly zero on $(n,\infty)$. Next, we impose Condition~\ref{cond:majorant}. This particularly implies that $\eta$ and $\delta$ are continuous random variables. We denote by $f_{(\eta,\delta)}$ the joint density function of $(\eta,\delta)$, by $f_{\eta|\delta}$ the conditional density function of $\eta$ given $\delta$, and by $f_\eta$ and $f_\delta$ the marginal density functions of $\eta$ and $\delta$, respectively. Finally, we suppose that the extended state process $\tilde{Z}$ is semi-Markovian. Because of Proposition~\ref{prop:linkmarkov}(a), the transition intensities $\mu_{jk}^{(0)}(\cdot) := \mu_{jk}^{\mathcal{F}}(\cdot,E)$, $j,k\in S$, $j \neq k$, from Condition~\ref{cond:majorant} may then be cast as
\begin{align*}
\mu_{jk}^{(0)}(t) &= \mu_{jk}^{(0)}(t,U_{t-}), \quad \quad \quad \quad j \in \{1,\ldots,\varsigma\}, k \in S, j \neq k, \\
\mu_{\text{p}\text{d}}^{(0)}(t) &= \mu_{\text{p}\text{d}}^{(0)}(t,U^{\text{h}}_{t-},U^{\text{p}}_{t-},H_{t-})
\end{align*} 
for $t\geq0$, where according to Proposition~\ref{prop:smoothness_fixed_again},
\begin{align*}
\mu_{jk}^{(0)}(t+,s) &= \lim_{h\downarrow 0}\frac{1}{h} P(Z_{t+h}=k \, | \, Z_t=j, U_t=s), \quad j \in \{1,\ldots,\varsigma\}, k \in S, j \neq k, \\
\mu_{\text{p}\text{d}}^{(0)}(t+,s,r,\ell) &= \lim_{h\downarrow 0} P(Z_{t+h}=\text{d} \, | \, Z_t = \text{p}, U^{\text{h}}_t = s, U^{\text{p}}_t = r, H_t = \ell)
\end{align*}
The remaining transition intensities are constantly zero.

The next results provide Feynman-Kac formulas that can serve as the starting point for the development of numerical schemes for the state-wise prospective reserves under full information, i.e.\ $V_j^{(0)}$, $j \in S$.
\begin{proposition}\label{prop:ex_F}
Suppose that the assumptions from the beginning of this section hold. If the function $W_{\emph{d}}^{(0)}(\cdot)$ is a bounded c\`{a}dl\`{a}g solution of
\begin{align}\label{semi_Markov_Thiele_in_J+2}
W_{\emph{d}}^{(0)}(\d t) = W_{\emph{d}}^{(0)}(t-)  \frac{\kappa(\d t)}{\kappa(t-)} - b_{\emph{d}}(t) \, \beta(\d t)
\end{align}
with terminal condition $W_{\emph{d}}^{(0)}(n) = 0$, and if the function $W^{(0)}_{\emph{p}}(\cdot,\cdot,\cdot,\cdot)$ is a bounded c\`{a}dl\`{a}g solution of
\begin{align}\label{semi_Markov_Thiele_in_J+1}\begin{split}
  W^{(0)}_{\emph{p}}(\d t, s,r,k) = & \, W^{(0)}_{\emph{p}}(t-,s,r,k) \frac{\kappa(\d t)}{\kappa(t-)} - b_{\emph{p}}(t) \, \beta(\d t) \\
  &  -  \big(b_{\emph{p}\emph{d}}(t) +  W^{(0)}_{\emph{d}}(t) - W^{(0)}_{\emph{p}}(t,s,r,k)\big) \, \mu_{\emph{p}\emph{d}}(t, t-s,t-r,k) \, \d t
\end{split}\end{align}
on $\{t>r>s\geq0, k\leq\varsigma\}$ with terminal conditions $W^{(0)}_{\emph{p}}(n,s,r,k)=0$ for $n \geq r \geq s \geq 0$, $k\leq\varsigma$, and if the functions $W^{(0)}_j(\cdot,\cdot)$, $j\in\{1,\ldots,\varsigma\}$, are bounded and c\`{a}dl\`{a}g solutions of
\begin{align}\label{semi_Markov_Thiele_in_i_correction}\begin{split}
   W^{(0)}_{j}(\d t,s) =&  \, W^{(0)}_j(t-,s) \frac{\kappa(\d t)}{\kappa(t-)} - b_j(t) \, \beta(\d t) \\
  &- \sum_{k \leq \varsigma: j \neq k} \big( b_{jk}(t) + W^{(0)}_{k}(t,t) - W^{(0)}_{j}(t,s)\big) \, \mu^{(0)}_{jk}(t,t-s) \,\d t\\
  & -  \big(b_{j\emph{p}}(t) + W^{(0)}_{\emph{p}}(t,s,t,j) - W^{(0)}_{j}(t,s)\big) \, \mu^{(0)}_{j\emph{p}}(t,t-s) \, \d t\\
  &-  \big( b_{j\emph{d}}(t) +  W^{(0)}_{\emph{d}}(t) - W^{(0)}_{j}(t,s)\big) \, \mu^{(0)}_{j\emph{d}}(t,t-s) \, \d t
\end{split}\end{align}
on $\{t>s\geq0,j\leq\varsigma\}$ with terminal conditions $W^{(0)}_j(n,s) = 0$ for $n \geq s \geq 0$, then $(W_j^{(0)}(t))_{t\geq0,j\in S}$ given by
\begin{align*}
W_j^{(0)}(t)
:=
\begin{cases}
W_j^{(0)}(t,t-I_j(t-)U_{t-}) & \text{ for } j\in\{1,\ldots,\varsigma\}, \\
W_{\emph{p}}^{(0)}(t,t-U_{t-}^{\emph{h}},t-U_{t-}^{\emph{p}},H_{t-}) & \text{ for } j = \emph{p}, \\
W_{\emph{d}}^{(0)}(t) & \text{ for } j = \emph{d},
\end{cases}
\quad \quad t\geq0,
\end{align*}
are state-wise prospective reserves in the sense of Definition~\ref{def:stat_counters} under full information $\mathcal{G}^{(0)}=\mathcal{F}$.
\end{proposition}
\begin{remark}
We may replace $I_j(t-)U_{t-}$, $U_{t-}^{\text{h}}$, $U_{t-}^{\text{p}}$, $H_{t-}$ by $U_t$, $U_t^{\text{h}}$, $U_t^{\text{p}}$, $H_t$ and still obtain condition (1) of Definition~\ref{def:stat_counters}, but if $I_j(t-)U_{t-}$ is replaced by $U_t$, then condition (2) of Definition~\ref{def:stat_counters} fails to hold.\demormk
\end{remark}
\begin{proof}
Note that  the right-continuity of the solutions of the differential equations allows us to uniquely expand the domains of the solutions to $t \geq s \geq 0$, $t\geq r > s \geq 0$, and $t \geq 0$. That means that $W^{(0)}_{j}(t,t)$ and $W^{(0)}_{\text{p}}(t,s,t,k)$ are indeed given by the solutions.

Note that $W^{(0)}_j = (W^{(0)}_j(t))_{t\geq0}$, $j\in S$, are $\mathcal{F}$-predictable bounded processes for which $t\mapsto \int_{(0,t]} I_j(s-) \, W^{(0)}_j(\d s)$ defines finite variation processes. So the result follows from Proposition~\ref{Prop:true_Thiele_equations} if we can show that~\eqref{eq:ThieleBSDE} is satisfied for every $j\in S$.

By multiplying~\eqref{semi_Markov_Thiele_in_J+2} with $I_{\text{d}}(t-)$, we obtain~\eqref{eq:ThieleBSDE} for $j = \text{d}$. By multiplying equation~\eqref{semi_Markov_Thiele_in_J+1} with $I_{\text{p}}(t-)$ and replacing $s$, $r$, $k$ by $t-U^{\text{h}}_{t-}$, $t-U^{\text{p}}_{t-}$, $H_{t-}$, we find that $W^{(0)}_{\text{p}}$ solves~\eqref{eq:ThieleBSDE} for $j=\text{p}$.

Let $\tau_i$, $i \in \mathbb{N}_0$, be the jump times of $Z$ with $\tau_0 := 0$. Multiplying equation~\eqref{semi_Markov_Thiele_in_i_correction} with $I_j(t-)\mathbf{1}_{\{\tau_i < t \leq \tau_{i+1}\}}$ and replacing $s$ by $\tau_i I_j(\tau_i)+\, t \,\mathbf{1}_{\{Z_{\tau_i}\neq j\}}$, we also obtain that $W^{(0)}_j(t, \tau_i I_j(\tau_i)+ t \,\mathbf{1}_{\{Z_{\tau_i}\neq j\}})$, $j\in \{1,\ldots, \varsigma\}$, is a solution of \eqref{eq:ThieleBSDE} on the interval $(\tau_i,\tau_{i+1}]$. This follows from
\begin{align*}
  &I_j(t-) W^{(0)}_{k}(t,\tau_i I_k(\tau_i)+\, t \,\mathbf{1}_{\{Z_{\tau_i}\neq k\}} ) = I_j(t-) W^{(0)}_{k}(t,t),\\
  & I_j(t-) W^{(0)}_{\text{p}}(t,\tau_i I_j(\tau_i)+ t \,\mathbf{1}_{\{Z_{\tau_i}\neq j\}},0,j) = I_j(t-) W^{(0)}_{\text{p}}(t,t- U^{\text{h}}_{t-} ,t-U_{t-}^{\text{p}},H_{t-})
\end{align*}
for all $t\in(\tau_i,\tau_{i+1}]$ and $j,k \in \{1, \ldots, \varsigma\}$, $j \neq k$. Summing over $i\in\mathbb{N}_0$ yields that the processes
\begin{align*}
  W^{(0)}_j(t, t- I_j(t-)U_{t-})= \sum_{i=0}^\infty  \mathbf{1}_{\{\tau_i < t \leq \tau_{i+1}\}} W^{(0)}_j(t, \tau_i I_j(\tau_i)+ t \,\mathbf{1}_{\{Z_{\tau_i}\neq j\}}), \quad j \in \{1,\ldots, \varsigma\},
\end{align*}
are solutions of~\eqref{eq:ThieleBSDE} for $j\in\{1,\ldots,\varsigma\}$ since
\begin{align*}
  I_j(t-) W^{(0)}_{k}(t,t- I_k(t-)U_{t-}) = I_j(t-)W^{(0)}_{k}(t,t)
\end{align*}
for all $t>0$ and $j,k \in \{1,\ldots,\varsigma\}$, $j \neq k$.
\end{proof}
The numerical schemes that can be developed based on Proposition~\ref{prop:ex_F} are significantly more complex than in the classic (semi-)Markovian case, see e.g.~\citet{BuchardtMollerSchmidt2015} and \citet{AdeChristiansen2017}. The sum at risks involve $W^{(0)}_{\text{p}}(t,s,t,j)$, which must be computed based on~\eqref{semi_Markov_Thiele_in_J+1} for all $0\leq s < t$.

For $i\in\{1,2\}$ we set $W_j^{(i)} := W_j^{(0)}$, $j\in\{1,\ldots,\varsigma,\text{d}\}$. The next results provide Feynman-Kac formulas under non-monotone information $\mathcal{G}^{(1)}$ and $\mathcal{G}^{(2)}$.

\begin{proposition}\label{prop:det_Thiele_eq_G1}
Suppose the assumptions from the beginning of this section hold. Suppose the function $W^{(1)}_{\emph{p}}(\cdot,\cdot)$ is a bounded and c\`{a}dl\`{a}g solution of
\begin{align}\label{eq:Thiele_J+1_G1}\begin{split}
W^{(1)}_{\emph{p}}(\d t,r) = \, &W^{(1)}_{\emph{p}}(t-,r) \frac{\kappa(\d t)}
{\kappa(t-)} - b_{\emph{p}}(t) \, \beta(\d t) - \big( b_{\emph{p}\emph{d}}(t)+W^{(1)}_{\emph{d}}(t)  - W^{(1)}_{\emph{p}}(t,r) \big) \, \mu^{(1)}_{\emph{p}\emph{d}}(t,r) \, \d t
\end{split}\end{align}
on $\{t > r > 0\}$ with terminal conditions $W^{(1)}_{\emph{p}}(n,r) = 0$ for $n\geq r \geq 0$, and where under the convention $0/0:=0$,
\begin{align}\label{formula:mortality_of_retirees_case1}
    \mu^{(1)}_{\emph{p}\emph{d}}(t,r) := \frac{ f_{\delta|\eta}(t|r)}{P( \delta \geq t \, | \, \eta=r)}, \quad t \geq r \geq 0.
\end{align}
Let $(W_{\emph{p}}^{(1)}(t))_{t\geq 0}$ be given by $W_{\emph{p}}^{(1)}(t) := W_{\emph{p}}^{(1)}(t,\eta)$ for $t\geq0$. Then $(W_j^{(1)}(t))_{t\geq0,j\in S}$ are state-wise prospective reserves in the sense of Definition~\ref{def:stat_counters} under non-monotone information $\mathcal{G}^{(1)}$.
\end{proposition}
\begin{proposition}\label{prop:det_Thiele_eq_G2}
Suppose the assumptions from the beginning of this section hold. If the function $W^{(2)}_{\emph{p}}(\cdot)$ is a bounded and c\`{a}dl\`{a}g solution of
\begin{align}\label{eq:Thiele_J+1_G2}\begin{split}
W^{(2)}_{\emph{p}}(\d t) = \, W^{(2)}_{\emph{p}}(t-) \frac{\kappa(\d t)}
{\kappa(t-)} - b_{\emph{p}}(t) \, \beta(\d t)
&- \big(b_{\emph{p}\emph{d}}(t)+W^{(2)}_{\emph{d}}(t)  - W^{(2)}_{\emph{p}}(t) \big) \, \mu^{(2)}_{\emph{p}\emph{d}}(t) \, \d t \\
  & + \big(W^{(1)}_{\emph{p}}(t,t)-W^{(2)}_{\emph{p}}(t)\big) \, \overline{\mu}^{(2)}_{\textrm{{\tiny $\bullet$}} \emph{p}}(t) \, \d t
\end{split}\end{align}
with terminal condition $W^{(2)}_{\emph{p}}(n) = 0$, where under the convention $0/0:=0$,
\begin{align}\label{formula:mortality_of_retirees_case2}
   \mu^{(2)}_{\emph{p}\emph{d}}(t) &:= \frac{\int_{0}^t f_{(\eta,\delta)}(s,t) \, \d s }{P(\eta < t \leq \delta)} ,\\
   \label{formula:mortality_of_retirees_adjust}
     \overline{\mu}^{(2)}_{\textrm{{\tiny $\bullet$}}\emph{p}}(t)  &:= \frac{f_{\eta}(t)}{P(\eta \leq t < \delta)},
\end{align}
then $(W_j^{(2)}(t))_{t\geq0,j\in S}$ are state-wise prospective reserves in the sense of Definition~\ref{def:stat_counters} under non-monotone information $\mathcal{G}^{(2)}$.
\end{proposition}
In order to reduce the computation time and simplify actuarial modeling and statistical estimation,  practitioners, when computing the prospective reserve for non-retirees based on $W^{(0)}_j$, $j \in \{1,\ldots,\varsigma\}$, often approximate $W^{(0)}_{\text{p}}(t,s,t,j)$ by a less complex quantity such as $W^{(1)}_{\text{p}}(t,t)$, which discards information concerning previous health records, or $W^{(2)}_{\text{p}}(t)$, which additionally discards information concerning the time of retirement. Replacing $W^{(0)}_{\text{p}}$ by $W^{(1)}_{\text{p}}$ or $W^{(2)}_{\text{p}}$ produces approximation errors on the individual level (and redistribution of wealth on the portfolio level for non-retirees).

Proposition~\ref{prop:det_Thiele_eq_G1} and Proposition~\ref{prop:det_Thiele_eq_G2} can be used to develop computational schemes for $W^{(1)}_{\text{p}}$ and $W^{(2)}_{\text{p}}$, respectively. Focusing on $W^{(2)}_{\text{p}}$, this involves the transition intensity $\mu^{(2)}_{\text{p}\text{d}}$, which by~\eqref{formula:mortality_of_retirees_case2} is the hazard rate corresponding to a classic mortality table for retirees. It also involves the adjustment term
\begin{align*}
\big(W^{(1)}_{\text{p}}(t,t)-W^{(2)}_{\text{p}}(t)\big) \, \overline{\mu}^{(2)}_{\textrm{{\tiny $\bullet$}}\text{p}}(t) \, \d t,
\end{align*}
where according to~\eqref{formula:mortality_of_retirees_adjust}, the term $ \overline{\mu}^{(2)}_{\textrm{{\tiny $\bullet$}}p}(t) \, \d t$ gives the infinitesimal probability of retirement having just occurred (at time $t$), conditionally on the insured presently being retired.

If the mortality does not depend on the time since retirement, i.e.\ if $\mu^{(1)}_{\text{p}\text{d}}(t,r) \equiv \mu^{(2)}_{\text{p}\text{d}}(t)$, we end up with the differential equations
\begin{align}\label{eq:practice_favorite}\begin{split}
W^{(i)}_{\text{p}}(\d t) = \, W^{(i)}_{\text{p}}(t-) \frac{\kappa(\d t)}
{\kappa(t-)} - b_{\text{p}}(t) \, \beta(\d t)
- \big(b_{\text{p}\text{d}}(t)+W^{(i)}_{\text{d}}(t)  - W^{(i)}_{\text{p}}(t) \big) \, \mu^{(2)}_{\text{p}\text{d}}(t) \, \d t
\end{split}\end{align}
for $i\in\{1,2\}$. Even though the mortality of retirees might depend on the time since retirement, practitioners often still utilize~\eqref{eq:practice_favorite} directly. This produces additional approximation errors on the individual level (and redistribution of wealth on the portfolio level for retirees and non-retirees).
\begin{proof}[Proof of Proposition~\ref{prop:det_Thiele_eq_G1}]
The assumptions guarantee that conditions (2) and (3) of Definition~\ref{def:stat_counters} are satisfied.  What remains to be shown is then condition (1): $I_{\text{p}}(t)W^{(1)}_{\text{p}}(t,\eta) = I_{\text{p}}(t)V^{(1)}(t)$ for all $t\geq0$. By applying integration by parts, we find that
\begin{align*}
  &\mathbf{1}_{\{\eta < t\}} \, \d \Big( I_{\text{p}}(t) v(t) W^{(1)}_{\text{p}}(t,\eta) \Big)\\
  &=  I_{\text{p}}(t-) \Big( v(t) W^{(1)}_{\text{p}}(\d t,\eta) -  v(t) W^{(1)}_{\text{p}}(t-,\eta) \frac{\kappa(\d t)}{\kappa(t-)}  \Big) - v(t) W^{(1)}_{\text{p}}(t,\eta) N_{\text{p}\text{d}}(\d t).
\end{align*}
Inserting~\eqref{eq:Thiele_J+1_G1} into the latter term leads to
  \begin{align*}
  &\mathbf{1}_{\{\eta < t\}} \, \d \Big( I_{\text{p}}(t) v(t) W^{(1)}_{\text{p}}(t,\eta) \Big)\\
  &=  - I_{\text{p}}(t-) v(t) B(\d t) - v(t) W_{\text{d}}^{(1)}(t) N_{\text{p}\text{d}}(\d t)    + v(t) r^{(1)}_{\text{p}\text{d}}(t)  M^{(1)}_{\text{p}\text{d}}(\d t),
\end{align*}
where $r^{(1)}_{\text{p}\text{d}}(t):=b_{\text{p}\text{d}}(t) + W^{(1)}_{\text{d}}(t)  - W^{(1)}_{\text{p}}(t,\eta)$, $t\geq0$, and $M^{(1)}_{\text{p}\text{d}}(\d t) := N_{\text{p}\text{d}}(\d t)-I_{\text{p}}(t-) \mu^{(1)}_{\text{p}\text{d}}(t, \eta) \, \d t$.
Thus, since $\{\eta < t\}\subseteq \{\eta < s\}$ for $s\geq t\geq 0$, it holds that
\begin{align*}
&\mathbf{1}_{\{\eta < t\}}  I_{\text{p}}(t) v(t) W^{(1)}_{\text{p}}(t,\eta) \\
&=
\E[\mathbf{1}_{\{\eta < t\}}  I_{\text{p}}(t) v(t) W^{(1)}_{\text{p}}(t,\eta)\, | \, \mathcal{G}^{(1)}_t] \\
&=
\mathbf{1}_{\{\eta < t\}} \E\bigg[v(t) \int_{(t,n]} I_{\text{p}}(s-)  \frac{\kappa (t)}{\kappa(s)} \, B(\d s) \, \bigg|  \,\mathcal{G}^{(1)}_t \bigg]+
\mathbf{1}_{\{\eta < t\}} \E\bigg[\int_{(t,n]}  v(s) W^{(1)}_{\text{d}}(s) \, N_{\text{p}\text{d}}(\d s) \, \bigg|  \,\mathcal{G}^{(1)}_t \bigg] \\
&\quad \, -
\mathbf{1}_{\{\eta < t\}} \E\bigg[ \int_{(t,n]} v(s) r^{(1)}_{\text{p}\text{d}}(s) \, M^{(1)}_{\text{p}\text{d}}(\d s) \, \bigg| \, \mathcal{G}^{(1)}_t \bigg] \\
&=
\mathbf{1}_{\{\eta < t\}} I_{\text{p}}(t) v(t) V^{(1)}(t)  -
\mathbf{1}_{\{\eta < t\}} I_{\text{p}}(t) \E\bigg[ \int_{(t,n]} v(s) r^{(1)}_{pd}(s) M^{(1)}_{pd}(\d s) \, \bigg| \, \mathcal{G}^{(1)}_t \bigg].
\end{align*}
The latter conditional expectation is zero since $\mathcal{G}^{(1)}_t$ equals the sub-filtration $\sigma(\mathbf{1}_{\{\eta \leq s\}}, \mathbf{1}_{\{\delta \leq s\}} : s \leq t)$ for $t > \eta$ and $v(s) r^{(1)}_{pd}(s) M^{(1)}_{pd}(\d s)$ is a martingale with respect to this sub-filtration. All in all, we conclude that since $v>0$,
\begin{align*}
\mathbf{1}_{\{\eta < t\}} I_{\text{p}}(t) W^{(1)}_{\text{p}}(t,\eta)
=
\mathbf{1}_{\{\eta < t\}} I_{\text{p}}(t) V^{(1)}(t)
\end{align*}
for all $t\geq0$. Recall that $V^{(1)}$ and $W^{(1)}_{\text{p}}$ c\`{a}dl\`{a}g sample paths. This allows us to replace $\mathbf{1}_{\{\eta < t\}} I_{\text{p}}(t)$ by $\mathbf{1}_{\{\eta \leq t\}} I_{\text{p}}(t) = I_{\text{p}}(t)$. Thus
\begin{align*}
I_{\text{p}}(t) W^{(1)}_{\text{p}}(t,\eta) = I_{\text{p}}(t) V^{(1)}(t)
\end{align*}
for all $t\geq0$ as desired.
\end{proof}

\begin{proof}[Proof of Proposition~\ref{prop:det_Thiele_eq_G2}]
The assumptions guarantee that conditions (2) and (3) of Definition~\ref{def:stat_counters} are satisfied.  What remains to be shown is then condition (1): $I_{\text{p}}(t)W^{(2)}_{\text{p}}(t) = I_{\text{p}}(t)V^{(2)}(t)$ for all $t\geq0$. Let $V^{(2)}_{\text{p}}$ be the state-wise prospective reserve from Theorem~\ref{cor:sdeY_G2}. So we can alternatively show $I_{\text{p}}(t)W^{(2)}_{\text{p}}(t) = I_{\text{p}}(t)V_{\text{p}}^{(2)}(t)$ for all $t\geq0$. Since $\eta$ and $\delta$ are continuous random variables, starting from~\eqref{remark:backward_sum_at_risk} it is possible to show that
\begin{align*}
\sum_{k \leq \varsigma}\int_E \overline{R}_{k\text{p}}^{(2)}(t,g) \, \overline{\lambda}_{k\text{p}}^{(2)}(\d t \times \d g)
=
\big(W_{\text{p}}^{(1)}(t,t) - V_{\text{p}}^{(2)}(t)\big) \, \overline{\mu}^{(2)}_{\textrm{{\tiny $\bullet$}}\text{p}}(t)\, \d t.
\end{align*}
Similarly we may show that
\begin{align*}
 \big(b_{\text{p}\text{d}}(t) + V_{\text{p}\text{d}}^{(2)}(t) - V_{\text{p}\text{p}}^{(2)}(t)\big) \lambda_{\text{p}\text{d}}^{(2)}(\d t)
 =
  \big(b_{\text{p}\text{d}}(t) + W_{\text{d}}^{(2)}(t) - V_{\text{p}}^{(2)}(t)\big) \, \mu_{\text{p}\text{d}}^{(2)}(t) \, \d t.
\end{align*}
Thus by applying integration by parts, inserting~\eqref{eq:Thiele_J+1_G2}, and applying Theorem~\ref{cor:sdeY_G2}, straightforward calculations yield
\begin{align*}
&\d \Big( I_{\text{p}}(t) v(t) W^{(2)}_{\text{p}}(t)-I_{\text{p}}(t)  v(t) V^{(2)}_{\text{p}}(t) \Big)^{\!2} \\
&= v(t)^2 I_{\text{p}}(t-) \big( W^{(2)}_{\text{p}}(t) - V_{\text{p}}^{(2)}(t)\big)^{2} \big(\mu^{(2)}_{\text{p}\text{d}}(t)-\overline{\mu}^{(2)}_{\textrm{{\tiny $\bullet$}} \text{p}}(t) \big) \,\d t\\
  & \quad -   v(t)^2\big( W^{(2)}_{\text{p}}(t) - V_{\text{p}}^{(2)}(t)\big)^{\!2} \big( N_{\text{p}\text{d}}(\d t)  -I_{\text{p}}(t-) \mu^{(2)}_{\text{p}\text{d}}(t) \, \d t \big)\\
  &\quad +   v(t)^2\big( W^{(2)}_{\text{p}}(t) - V_{\text{p}}^{(2)}(t)\big)^{\!2} \bigg(\!-I_{\text{p}}(t)\overline{\mu}^{(2)}_{\textrm{{\tiny $\bullet$}} \text{p}}(t)  \, \d t + \sum_{k\neq\varsigma} N_{k\text{p}}(\d t)\bigg).
\end{align*}
The last two terms are infinitesimal forward and backward martingales, respectively. It is possible to show that they equal zero in expectation. Then
\begin{align*}
  &v(t)^2 P(Z_t = \text{p})\big( W^{(2)}_{\text{p}}(t) - V_{\text{p}}^{(2)}(t)\big)^{\!2} \\
  &= \E\Big[  I_{\text{p}}(t) v(t)^2 \big( W^{(2)}_{\text{p}}(t) - V_{\text{p}}^{(2)}(t)\big)^{\!2} \, \Big]\\
  &= - \E \bigg[\int_{(t,n]} v(s)^2 I_{\text{p}}(s-) \big( W^{(2)}_{\text{p}}(t) - V_{\text{p}}^{(2)}(t)\big)^{\!2} \big(\mu^{(2)}_{\text{p}\text{d}}(s)-\overline{\mu}^{(2)}_{\textrm{{\tiny $\bullet$}} \text{p}}(s)\big) \,\d s\bigg] \\
  &= - \int_{(t,n]} v(s)^2 P(Z_s = p) \big( W^{(2)}_{\text{p}}(t) - V_{\text{p}}^{(2)}(t)\big)^{\!2} \big(\mu^{(2)}_{\text{p}\text{d}}(s)-\overline{\mu}^{(2)}_{\textrm{{\tiny $\bullet$}} \text{p}}(s)\big) \, \d s.
\end{align*}
This means that the function
\begin{align*}
g(t) :=  v(t)^2 P(Z_t = \text{p}) \big( W^{(2)}_{\text{p}}(t) - V_{\text{p}}^{(2)}(t)\big)^{\!2}
\end{align*}
satisfies the integral equation
\begin{align*}
g(t) = -\int_{(t,n]} g(s) \big(\mu^{(2)}_{\text{p}\text{d}}(s)-\overline{\mu}^{(2)}_{\textrm{{\tiny $\bullet$}} \text{p}}(s)\big) \, \d s
\end{align*}
for all $t\in[0,n]$ under the convention $(n,n] = \emptyset$. Note that
\begin{align*}
\lvert g(t) \rvert \leq
\int_{(t,n]} \lvert g(s) \rvert \, \big\lvert \mu^{(2)}_{\text{p}\text{d}}(s)-\overline{\mu}^{(2)}_{\textrm{{\tiny $\bullet$}} \text{p}}(s) \big\rvert \, \d s
\end{align*}
for all $t \in [0,n]$. According to the (path-wise) backward Gr\"{o}nwall inequality \citep[see Lemma 4.7 in][]{CohenElliot2012}, $g(t) = 0$ for all $t\in[0,n]$. Since $v>0$, for each $t\geq 0$ it then holds that $I_{\text{p}}(t) W^{(2)}_{\text{p}}(t)= I_{\text{p}}(t) V_{\text{p}}^{(2)}(t)$ as desired.
\end{proof}

\section*{Acknowledgments and declarations of interest}

Christian Furrer's research has partly been funded by the Innovation Fund Denmark (IFD) under File No.\ 7038-00007B. The authors declare no competing interests.


\begin{thebibliography}{29}
\providecommand{\natexlab}[1]{#1}
\providecommand{\url}[1]{\texttt{#1}}
\expandafter\ifx\csname urlstyle\endcsname\relax
  \providecommand{\doi}[1]{doi: #1}\else
  \providecommand{\doi}{doi: \begingroup \urlstyle{rm}\Url}\fi

\bibitem[Aalen(1978)]{Aalen1978}
O.~Aalen (1978).
\newblock {Nonparametric Inference for a Family of Counting Processes}.
\newblock \emph{Annals of Statistics}, 6\penalty0 (4):\penalty0 701--726.
\newblock \doi{10.1214/aos/1176344247}.

\bibitem[Ad{\'e}kambi and Christiansen(2017)]{AdeChristiansen2017}
F.~Ad{\'e}kambi and M.C. Christiansen (2017).
\newblock Integral and differential equations for the moments of multistate
  models in health insurance.
\newblock \emph{Scandinavian Actuarial Journal}, 2017:\penalty0 29--50.
\newblock \doi{10.1080/03461238.2015.1058854}.

\bibitem[Bladt et~al.(2020)Bladt, Asmussen, and Steffensen]{Bladt2020}
M.~Bladt, S.~Asmussen, and M.~Steffensen (2020).
\newblock Matrix representations of life insurance payments.
\newblock \emph{European Actuarial Journal}.
\newblock \doi{10.1007/s13385-019-00222-0}.

\bibitem[Buchardt et~al.(2015)Buchardt, M{\o}ller, and
  Schmidt]{BuchardtMollerSchmidt2015}
K.~Buchardt, T.~M{\o}ller, and K.B. Schmidt (2015).
\newblock Cash flows and policyholder behaviour in the semi-{M}arkov life
  insurance setup.
\newblock \emph{Scandinavian Actuarial Journal}, 2015\penalty0 (8):\penalty0
  660--688.
\newblock \doi{10.1080/03461238.2013.879919}.

\bibitem[Christiansen(2008{\natexlab{a}})]{Christiansen2008a}
M.C. Christiansen (2008{\natexlab{a}}).
\newblock A sensitivity analysis concept for life insurance with respect to a
  valuation basis of infinite dimension.
\newblock \emph{Insurance: Mathematics and Economics}, 42\penalty0
  (2):\penalty0 680--690.
\newblock \doi{10.1016/j.insmatheco.2007.07.005}.

\bibitem[Christiansen(2008{\natexlab{b}})]{Christiansen2008b}
M.C. Christiansen (2008{\natexlab{b}}).
\newblock A sensitivity analysis of typical life insurance contracts with
  respect to the technical basis.
\newblock \emph{Insurance: Mathematics and Economics}, 42\penalty0
  (2):\penalty0 787--796.
\newblock \doi{10.1016/j.insmatheco.2007.08.005}.

\bibitem[Christiansen(2021)]{Christiansen2020}
M.C. Christiansen (2021).
\newblock A martingale concept for non-monotone information in a jump process
  framework.
\newblock arXiv: \href{https://arxiv.org/abs/1811.00952}{1811.00952}.

\bibitem[Christiansen and Djehiche(2020)]{ChristiansenDjehiche2020}
M.C. Christiansen and B.~Djehiche (2020).
\newblock Nonlinear reserving and multiple contract modifications in life
  insurance.
\newblock \emph{Insurance: Mathematics and Economics}, 93:\penalty0 187--195.
\newblock \doi{10.1016/j.insmatheco.2020.05.004}.

\bibitem[Christiansen and Steffensen(2013)]{ChristiansenSteffensen2013}
M.C. Christiansen and M.~Steffensen (2013).
\newblock {SAFE-SIDE SCENARIOS FOR FINANCIAL AND BIOMETRICAL RISK}.
\newblock \emph{ASTIN Bulletin}, 43\penalty0 (3):\penalty0 323–357.
\newblock \doi{10.1017/asb.2013.16}.

\bibitem[Cohen and Elliott(2012)]{CohenElliot2012}
S.N. Cohen and R.J. Elliott (2012).
\newblock Existence, uniqueness and comparisons for {BSDE}s in general spaces.
\newblock \emph{The Annals of Probability}, 40\penalty0 (5):\penalty0
  2264--2297.
\newblock \doi{10.1214/11-AOP679}.

\bibitem[Djehiche and L{\"o}fdahl(2016)]{DjehicheLofdahl2016}
B.~Djehiche and B.~L{\"o}fdahl (2016).
\newblock Nonlinear reserving in life insurance: Aggregation and mean-field
  approximation.
\newblock \emph{Insurance: Mathematics and Economics}, 69\penalty0
  (3):\penalty0 1--13.
\newblock \doi{10.1016/j.insmatheco.2016.04.002}.

\bibitem[{European Parliament} and {Council of the European Union}(2016)]{GDPR}
{European Parliament} and {Council of the European Union} (2016).
\newblock Regulation ({EU}) 2016/679 of the {E}uropean {P}arliament and of the
  {C}ouncil of 27 {A}pril 2016 on the protection of natural persons with regard
  to the processing of personal data and on the free movement of such data, and
  repealing {D}irective 95/46/{EC} ({G}eneral {D}ata {P}rotection
  {R}egulation).
\newblock
  \url{https://eur-lex.europa.eu/legal-content/EN/TXT/PDF/?uri=CELEX:32016R0679}.

\bibitem[Helwich(2008)]{Helwich2008}
M.~Helwich (2008).
\newblock \emph{Durational effects and non-smooth semi-{M}arkov models in life
  insurance}.
\newblock PhD thesis, University of Rostock.

\bibitem[Hoem(1969)]{Hoem1969}
J.M. Hoem (1969).
\newblock {Markov Chain Models in Life Insurance}.
\newblock \emph{Bl{\"a}tter der DGVFM}, 9:\penalty0 91--107.
\newblock \doi{10.1007/BF02810082}.

\bibitem[Jacobsen(2006)]{Jacobsen2006}
M.~Jacobsen (2006).
\newblock \emph{Point process theory and applications: Marked point and
  piecewise deterministic processes}.
\newblock Probability and its Applications. Birkh{\"a}user.
\newblock \doi{10.1007/0-8176-4463-6}.

\bibitem[Jarner and M{\o}ller(2015)]{JarnerMoller2015}
S.F. Jarner and T.~M{\o}ller (2015).
\newblock A partial internal model for longevity risk.
\newblock \emph{Scandinavian Actuarial Journal}, 2015\penalty0 (4):\penalty0
  352--382.
\newblock \doi{10.1080/03461238.2013.836561}.

\bibitem[Kalashnikov and Norberg(2003)]{KalashnikovNorberg2003}
V.~Kalashnikov and R.~Norberg (2003).
\newblock {On the Sensitivity of Premiums and Reserves to Changes in Valuation
  Elements}.
\newblock \emph{Scandinavian Actuarial Journal}, 2003\penalty0 (3):\penalty0
  238--256.
\newblock \doi{10.1080/03461230110106408}.

\bibitem[Linnemann(1993)]{Linnemann1993}
P.~Linnemann (1993).
\newblock On the application of {T}hiele's differential equation in life
  insurance.
\newblock \emph{Insurance: Mathematics and Economics}, 13\penalty0
  (1):\penalty0 63--74.
\newblock \doi{10.1016/0167-6687(93)90536-X}.

\bibitem[Milbrodt and Stracke(1997)]{MilbrodtStracke1997}
H.~Milbrodt and A.~Stracke (1997).
\newblock Markov models and {T}hiele's integral equations for the prospective
  reserve.
\newblock \emph{Insurance: Mathematics and Economics}, 19\penalty0
  (3):\penalty0 187--235.
\newblock \doi{10.1016/S0167-6687(97)00020-6}.

\bibitem[M{\o}ller(1993)]{Moller1993}
C.M. M{\o}ller (1993).
\newblock A stochastic version of {T}hiele's differential equation.
\newblock \emph{Scandinavian Actuarial Journal}, 1993:\penalty0 1--16.
\newblock \doi{10.1080/03461238.1993.10413910}.

\bibitem[Norberg(1991)]{Norberg1991}
R.~Norberg (1991).
\newblock {Reserves in Life and Pension Insurance}.
\newblock \emph{Scandinavian Actuarial Journal}, 1991:\penalty0 3--24.
\newblock \doi{10.1080/03461238.1991.10557357}.

\bibitem[Norberg(1992)]{Norberg1992}
R.~Norberg (1992).
\newblock Hattendorff's theorem and {T}hiele's differential equation
  generalized.
\newblock \emph{Scandinavian Actuarial Journal}, 1992:\penalty0 2--14.
\newblock \doi{10.1080/03461238.1992.10413894}.

\bibitem[Norberg(1996)]{Norberg1996}
R.~Norberg (1996).
\newblock {Addendum to Hattendorff's Theorem and Thiele's Differential Equation
  Generalized, SAJ 1992, 2–14}.
\newblock \emph{Scandinavian Actuarial Journal}, 1996:\penalty0 50--53.
\newblock \doi{10.1080/03461238.1996.10413962}.

\bibitem[Norberg(1999)]{Norberg1999}
R.~Norberg (1999).
\newblock A theory of bonus in life insurance.
\newblock \emph{Finance and Stochastics}, 3:\penalty0 373--390.
\newblock \doi{10.1007/s007800050067}.

\bibitem[Norberg(2001)]{Norberg2001}
R.~Norberg (2001).
\newblock {On Bonus and Bonus Prognoses in Life Insurance}.
\newblock \emph{Scandinavian Actuarial Journal}, 2001\penalty0 (2):\penalty0
  126--147.
\newblock \doi{10.1080/03461230152592773}.

\bibitem[Ramlau-Hansen(1988{\natexlab{a}})]{RamlauHansen1988a}
H.~Ramlau-Hansen (1988{\natexlab{a}}).
\newblock {Hattendorff's Theorem: A Markov chain and counting process
  approach}.
\newblock \emph{Scandinavian Actuarial Journal}, 1988\penalty0 (1-3):\penalty0
  143--156.
\newblock \doi{10.1080/03461238.1988.10413845}.

\bibitem[Ramlau-Hansen(1988{\natexlab{b}})]{RamlauHansen1988b}
H.~Ramlau-Hansen (1988{\natexlab{b}}).
\newblock The emergence of profit in life insurance.
\newblock \emph{Insurance: Mathematics and Economics}, 7\penalty0 (4):\penalty0
  225--236.
\newblock \doi{10.1016/0167-6687(88)90080-7}.

\bibitem[Serfozo(1971)]{Serfozo1971}
R.F. Serfozo (1971).
\newblock {Functions of Semi-Markov Processes}.
\newblock \emph{SIAM Journal on Applied Mathematics}, 20\penalty0 (3):\penalty0
  530--535.
\newblock \doi{10.1137/0120055}.

\bibitem[Steffensen(2000)]{Steffensen2000}
M.~Steffensen (2000).
\newblock A no arbitrage approach to {T}hiele's differential equation.
\newblock \emph{Insurance: Mathematics and Economics}, 27\penalty0
  (2):\penalty0 201--214.
\newblock \doi{10.1016/S0167-6687(00)00048-2}.

\end{thebibliography}

\end{document}